\documentstyle{amsart}

\def\ra{\rightarrow}

\def\Z{\mbox{$\Bbb Z$}}
\def\N{\mbox{$\Bbb N$}}

\def\ph{\varphi}

\def\and{\hbox{ and }}
\def\s#1{\mbox{$(-1)^{#1}$}}
\def\e#1{\mbox{$|#1|$}}
\def\ip#1#2{\left<#1,#2\right>}
\def\br#1#2{\left\{#1,#2\right\}}
\def\brf{\left\{\cdot,\cdot\right\}}
\def\brt#1#2{\left\{#1,#2\right\}}
\def\brtf{\left\{\cdot,\cdot\right\}}

\def\ipf{\left<\cdot,\cdot\right>}
\def\blacksquare{\quad \vrule height 6pt width 6pt}
\def\inv{^{-1}}

\newtheorem{thm}{Theorem}

\newtheorem{lma}{Lemma}
\def\ie{\hbox{\it i.e.}}
\def\ainf{\hbox{$A_\infty$}}
\def\linf{\mbox{$L_\infty$}}
\def\zt{\hbox{$\Z_2$}}
\def\ztz{\hbox{$\Z_2\times\Z$}}
\def\hd{\hbox{{$\hat d$}}}
\def\hm{\hbox{$\hat m$}}
\def\hl{\hbox{$\hat l$}}
\def\hmu{\hbox{$\hat \mu$}}

\def\tf{\hbox{$\tilde f$}}
\def\tl{\hbox{$\tilde l$}}
\def\tph{\hbox{$\tilde \ph$}}

\def\bdl{\hbox{$\bar\delta$}}
\def\bmu{\hbox{$\bar\mu$}}
\def\brd{\hbox{$\bar d$}}
\def\bm{\hbox{$\bar m$}}

\def\hdl{\hbox{$\hat \delta$}}

\def\CW{C(W)}
\def\CV{C(V)}
\def\hom{\mbox{\rm Hom}}

\def\sb#1{[#1]}
\def\tens{\bigotimes}

\def\tns{\otimes}
\def\mtns{\tns\cdots\tns}
\def\modot{\odot\cdots\odot}

\def\mwedge{\wedge\cdots\wedge}
\def\modot{\odot\cdots\odot}
\def\mplus{+\cdots+}
\def\mcom{,\cdots,}
\def\iso{\kern.35em{\raise3pt\hbox{$\sim$}\kern-1.1em\to}
         \kern.3em}
\def\im{\hbox{im}}

\def\k{\mbox{\bf k}}

\def\V{\mbox{V}}

\def\bid{\operatorname{bid}}
\def\bd#1{\bid(#1)}
\def\shf{\operatorname{Sh}}
\def\sh#1#2{\shf(#1,#2)}

\def\deg{\operatorname{deg}}

\def\coder{\operatorname{Coder}}

\def\s#1{\mbox{$(-1)^{#1}$}}
\def\e#1{\mbox{$|#1|$}}
\def\ip#1#2{\left<#1,#2\right>}
\def\blacksquare{\quad \vrule height 6pt width 6pt}

\begin{document}
\nocite{mar,sta3,getz2,lod,conn,kast,seib,hoch,umb,aksz}
\author{Michael Penkava}
\thanks{Partially supported by NSF grant ?? The author would also like
to thank the University of Washington, Seattle for hosting him this
summer.}
\address{University of California\\
Davis, CA 95616}
\email{michae@@math.ucdavis.edu}
\title{\linf\ Algebras and their Cohomology}
\maketitle
\begin{abstract}
An associative multiplication structure is nothing but a special type
of odd codifferential on the tensor coalgebra of a vector space, and
an $A_\infty$ algebra is simply a more general type of codifferential.
Hochschild cohomology classifies the deformations of an associative
algebra into an $A_\infty$ algebra, and cyclic cohomology in the presence
of an invariant inner product classifies the deformations of the associative
algebra into an $A_\infty$ algebra preserving the inner product.

Similarly, a Lie algebra or superalgebra structure is simply a special
case of an odd codifferential on the exterior coalgebra of a vector space,
and a $L_\infty$ algebra is given by an arbitrary codifferential.
We define ordinary and cyclic cohomology of $L_\infty$ algebras.
Lie algebra cohomology
classifies deformations of a Lie algebra into an $L_\infty$ algebra.
and cyclic
cohomology classifies the deformations of the Lie algebra into an
$L_\infty$ algebra preserving an invariant inner product.
In the case of Lie algebras, the exterior coalgebra is dual to the
symmetric coalgebra of the parity reversion of the vector space, and
by use of the parity reversion one obtains a realization of the differential
in terms of a structure of a differential graded Lie algebra on the space of
cochains.
\end{abstract}
\section{Introduction}
In a joint paper with Albert Schwarz \cite{ps2}, we gave definitions
of Hochschild cohomology and cyclic cohomology
of an \ainf\ algebra, and showed that
these cohomology
theories classified the infinitesimal deformations of the
\ainf\ structure and those deformations preserving an invariant inner
product. Then we showed that the Hochschild cohomology of an associative
algebra classifies the deformations of the algebra into an \ainf\
algebra, and the cyclic cohomology of an algebra with an invariant inner
product classifies the deformations of the algebra into an \ainf\ algebra
preserving  the inner product.
We then applied these results to show that
cyclic cocycles of an associative algebra determine homology cycles
in the complex of metric ribbon graphs.

The original purpose
of this paper was to apply the same constructions to \linf\ algebras,
and use the results to obtain homology cycles in another graph complex,
that of ordinary metric graphs.
However, in the preparation of the paper, I found that many of the
notions which are needed to explain the results are not easy to find
in the literature. Thus, to make the article more self contained, I
decided to include definitions of cyclic cohomology of Lie algebras,
cohomology of \zt-graded algebras, and coderivations of the tensor,
exterior and symmetric coalgebras.
It also became apparent that
treatment of the cohomology of
\ainf\ algebras from the perspective of codifferentials
on the tensor coalgebra was useful; we avoided this description in
the joint article for simplicity, but we missed some important ideas
because of this fact.
{}From the coalgebra point of view, the
corresponding theory of \linf\ algebras is seen to be closely analogous
to the theory of \ainf\ algebras, so the formulation and proofs of the
results about the \linf\ case are given near the end of the text.
The ideas in this text lead immediately to a simple formulation of
the cycle in the complex of metric ribbon graphs  associated to an
\ainf\ algebra with an invariant inner product. This same method
can be applied to show that \linf\ algebras with an invariant inner
product give rise to a cycle in the homology of
the complex of metric ordinary graphs.  These results will
be presented in a paper to follow this.

The notion of an \ainf\ algebra,
also called a strongly homotopy associative algebra,
was introduced by J. Stasheff in \cite{sta1,sta2}, and is a
generalization of an associative algebra.  From a certain point of
view, an associative algebra is simply a special case of a
codifferential on the tensor coalgebra of a vector space.
An \ainf\ algebra is given by taking an arbitrary coderivation;
in particular
associative algebras  and differential graded associative algebras are
examples of \ainf\ algebras.

\linf\ algebras, also called strongly homotopy Lie algebras,
first appeared in \cite{ss}, and are generalizations of Lie algebras.
A Lie algebra can be viewed as simply a special case of a codifferential
on the exterior coalgebra of a vector space, and \linf\ algebras are
simply arbitrary codifferentials on this coalgebra.

A bracket structure was introduced on
the space of cochains of an associative algebra by M. Gerstenhaber
in \cite{gers}.
{}From the coagebra point of view, the Gerstenhaber bracket
turns out to be simply the bracket of coderivations.
The space of cochains of a Lie algebra with coefficients in the
adjoint representation
also has a natural bracket,
which is the bracket of cochains. In this paper, we generalize these
brackets to \ainf\ and \linf\ algebras,
and define a bracket for cyclic
cohomology when there is an invariant inner product.

In our considerations, we shall be interested in \zt-graded graded
spaces, but we should point out that all the results hold in the
\Z-graded case as well. \ainf\ algebras were first defined as
\Z-graded objects, but for the applications we have in mind, the
\zt-grading is more appropriate, and the generalization of the
results here to the \Z-graded case is straightforward.
We shall find it necessary to
consider the parity reversion of a \zt-graded space. This is the
same space with the parity of elements reversed.
(In the \Z-graded case, the corresponding notion is that of suspension.)
There is
a natural isomorphism between the tensor coalgebra of a \zt-graded space
and the tensor coalgebra of its parity reversion. But in the case of
the exterior coalgebra, the isomorphism is to the {\em symmetric}
coalgebra of the parity reversion, a subtle point that
is not clarified very well in the literature.

A notion that will play a crucial role in what follows is that of
a grading group. An abelian group $G$ is said to be a grading group
if it possesses a symmetric \zt-valued bilinear form $\ipf$. Any
abelian group with a subgroup of index 2 possesses a natural
grading form, but this is not always the form which we shall need
to consider. An element $g$ of $G$ is called odd if $\ip gg=1$.
A grading group with a nontrivial
inner product is called good if $\ip gh=1$, whenever $g$ and $h$ are
both odd. Groups equipped with the natural inner product induced by
a subgroup of index 2 are good, and these include both \zt\ and \Z.
If $G$ and $H$ are grading groups, then $G\times H$ has an induced
inner product, given by $\ip{(g,g')}{(h,h')}=\ip g{g'}+\ip h{h'}$.
But $G\times H$ is never good when $G$ and $H$ are good. Now, if
$V$ is a $G$ graded vector space, then one can define the symmetric
and exterior algebras of the tensor algebra $T(V)$. The symmetric
algebra is $G$-graded commutative, but the exterior algebra is not.
On the other hand, the tensor, symmetric and exterior algebras also
are graded by $G\times\Z$, and with respect to the induced inner
product on $G\times\Z$, the exterior algebra is graded commutative.
Thus the grading group associated to the exterior algebra is not good.
The consequences of this fact play an important role in the theory
of \ainf\  and \linf\ algebras.

The organization of the paper is as follows. In section \ref{sect 1}
we recall the definition of the coboundary operator for a Lie algebra,
give a definition of cyclic cochain, define cyclic cohomology for a
Lie algebra, and relate it to deformations preserving an invariant
inner product. Section \ref{sect 2} extends these notions to the
case of a \zt-graded Lie algebra. The main purpose of the first two
sections is to present the formulas for comparison to the later
generalized constructions. In section \ref{sect 3} we give definitions
of the exterior and symmetric algebras, and fix conventions we
use for parity and bidegree, and the inner products on the grading
groups. Section \ref{sect 4} explains the tensor, exterior and
symmetric coalgebra structures, and their natural grading groups.

Section \ref{sect 5} discusses coderivations
and codifferentials of these coalgebras, and the
dependence on whether the grading group is \zt\ or \ztz. Section
\ref{sect 6} examines the duality between the tensor coalgebra
of a space with \ztz-grading, and the tensor coalgebra of the
parity reversion of the space with the \zt-grading. The notion
of \ainf\ algebra is introduced, the Gerstenhaber bracket is defined,
and it is used to define
cohomology of an \ainf\ algebra.
Invariant inner products are
defined, and the bracket of cyclic cochains is used to define
the coboundary of
cyclic cohomology. The bracket structure
on cyclic cochains depends on the inner product, but as usual for
cyclic cohomology, the  coboundary does not depend on it. Finally,
in section \ref{sect 7} the duality between the exterior coalgebra
of a space with the \ztz-grading, and the symmetric coalgebra of
the parity reversion of the space with the \zt-grading is used
to give a definition of an \linf\ algebra. The results of
the previous section are extended to the \linf\ algebra case.
\section{Cohomology of Lie Algebras}\label{sect 1}
\maketitle
In this section we recall the definition of the cohomology of a Lie algebra
with coefficients in a module, and relate the cohomology of the Lie
algebra with coefficients in the adjoint representation to
the
theory of deformations of the Lie algebra structure. Then we discuss
the theory of deformations of a Lie algebra preserving an invariant
inner product.

Let us suppose that $V$ is a Lie algebra over a field \k, with
bracket $\left[\cdot,\cdot\right]$.
The antisymmetry of the bracket means that the bracket is a linear map
from $\bigwedge^2\V\ra \V$.
Let $M$ be a $\V$ module, and
let $C^n(\V,M)=\hom(\bigwedge^n \V,M)$ be the space
of antisymmetric $n$-multilinear functions on $\V$ with values in $M$,
which we will call the module of cochains of degree $n$ on
$\V$ with values in $M$. The Lie algebra coboundary operator
$d:C^n(\V,M)\ra C^{n+1}(\V,M)$ is defined by
\begin{multline}
d\ph(v_1\mcom  v_{n+1})=\\
\sum_{1\le i<j\le n+1}\s{i+j-1}\ph([v_i,v_j],v_1\mcom
\hat v_i\mcom\hat v_j\mcom v_{n+1})\\
+
\sum_{1\le i\le n+1}\s iv_i\cdot\ph(v_1\mcom \hat v_i\mcom v_{n+1})
\end{multline}
Then
\begin{equation}
H^n(\V,M)=\ker(d:C^n(\V,M)\ra C^{n+1}(\V,M))/\im(d:C^{n-1}(\V,M)\ra C^n(\V,M))
\end{equation}
is the Lie algebra cohomology of $\V$ with coefficients in $M$. When
$M=\V$, and the action of $\V$ on itself is the adjoint action, then
we denote $C^n(\V,\V)$ as simply $C^n(\V)$, and similarly $H^n(\V,\V)$ is
denoted by $H^n(\V)$.

The connection between the cohomology of the Lie algebra and
(infinitesimal) deformations of $\V$ is given by $H^2(\V)$. If
we denote the bracket in $\V$ by $l$, and an infinitesimally
deformed product by $l_t=l+t\ph$, with $t^2=0$, then the map
$\ph:\bigwedge^2 \V\ra \V$ is a cocycle, and the trivial deformations
are coboundaries. For by the Jacobi identity we have
\begin{equation}
l_t(v_1,l_t(v_2,v_3))=l_t(l_t(v_1,v_2),v_3)+l_t(v_2,l_t(v_1,v_3))
\end{equation}
Expanding the expression above we determine that
\begin{multline}
[v_1,\ph(v_2,v_3)]+\ph(v_1,[v_2,v_3])=
[\ph(v_1,v_2),v_3]+\ph([v_1,v_2],v_3)\\
+[v_2,\ph(v_1,v_3)]+\ph(v_2,[v_1,v_3]),
\end{multline}
which can be expressed as the condition $d\ph=0$.

On the other hand, the notion of a trivial deformation is that $\V$ with
the new bracket is isomorphic to the original bracket structure. This
means that there is a linear bijection $\rho_t:A\ra A$ such that
$l_t(\rho_t(v_1),\rho_t(v_2))=\rho_t([v_1,v_2])$. We can express
$\rho_t=I+t\lambda$, where $\lambda:\V\ra \V$ is a linear map. Then
\begin{equation}
l_t(v_1,v_2)=l(v_1,v_2)+t(d\lambda)(v_1,v_2).
\end{equation}
Thus the trivial deformations are precisely those which are given by
coboundaries.

Next, consider an invariant
inner product on $\V$, by
which we mean a non-degenerate
symmetric bilinear form $\ip..:\V\tens \V\ra\k$ which satisfies
\begin{equation}
\ip{[v_1,v_2]}{v_3}=\ip{v_1}{[v_2,v_3]}.
\end{equation}
Note that for an invariant inner
product, the tensor $\tl$ given by
\begin{equation}
\tl(v_1,v_2,v_3)=\ip{[v_1,v_2]}{v_3}
\end{equation}
is also antisymmetric, so that $\tl\in\hom(\bigwedge^3 \V,\k)$. We also note
that $$\tl(v_1,v_2,v_3)=\tl(v_3,v_1,v_2),$$ so that $\tl$ is invariant
under cyclic permutations of $v_1$, $v_2$, $v_3$. We are interested in
deformations of $A$ which preserve this inner product, and these
deformations are given by $H^3(\V,\k)$, the cohomology of $\V$ with
trivial coefficients. To see this connection, we first define an
element $\ph\in C^n(\V)$ to be cyclic with respect to the inner product
if
\begin{equation}
\ip{\ph(v_1\mcom v_n)}{v_{n+1}}=\s{n}\ip{v_1}{\ph(v_2\mcom v_{n+1})}
\end{equation}
Then it is easy to see that $\ph$ is cyclic iff the map
$\tph:\V^n\ra \V$ given by
\begin{equation}
\tph(v_1\mcom v_{n+1})=
\ip{\ph(v_1\mcom v_n)}{v_{n+1}}
\end{equation}
 is antisymmetric, in other words,
$\tph\in\hom(\bigwedge^{n+1}\V,\k)$.
The term cyclic here is used to express the fact that $\tph$ is cyclic
in the sense that
\begin{equation}
\tph(v_1\mcom v_{n+1})=\s{n}\tph(v_{n+1},v_1\mcom v_n),
\end{equation}
which holds for any antisymmetric form.

Since the inner product is non-degenerate, the map
$\ph\mapsto\tilde\ph$ is an isomorphism between
the subspace $CC^n(\V)$ consisting of cyclic cochains,
and $C^{n+1}(\V,\k)$. We shall see that $d\ph$ is cyclic when $\ph$ is,
so that we can define the cyclic cohomology of the Lie algebra to be
\begin{equation}
HC^n(\V)=\ker(d:CC^n(\V)\ra CC^{n+1}(\V))/\im(d:CC^{n-1}(\V)\ra CC^n(\V)),
\end{equation}
We shall also show
that the isomorphism between $CC^n(\V)$ and
$C^{n+1}(\V,\k)$ commutes with the coboundary operator, so
that the cohomology of the complex of cyclic cochains coincides with
the cohomology of the Lie algebra with trivial coefficients, with
degree shifted by 1, \ie, $HC^n(\V)\cong H^{n+1}(\V,\k)$.
Thus, unlike the case of associative algebras,
cyclic cohomology does not lead to anything new.

To see these facts, note that if $\ph$ is cyclic, then
\begin{multline}
\ip{[v_i,\ph(v_1\mcom \hat v_i\mcom v_{n+1})]}{v_{n+2}}=
-\ip{[\ph(v_1\mcom \hat v_i\mcom v_{n+1}),v_i]}{v_{n+2}}\\
=-\ip{\ph(v_1\mcom \hat v_i\mcom v_{n+1})}{[v_i,v_{n+2}]}=
-\tph(v_1\mcom \hat v_i\mcom v_{n+1},[v_i,v_{n+2}])\\
=\s{n+1}\tph([v_i,v_{n+2}],v_1\mcom \hat v_i\mcom v_{n+1}).
\end{multline}
Thus
\begin{multline}
\widetilde{d\ph}(v_1\mcom v_{n+2})=
\ip{d\ph(v_1\mcom v_{n+1})}{v_{n+2}}=\\
\sum_{1\le i<j\le n+1}\s{i+j-1}\ip{\ph([v_i,v_j],v_1\mcom
\hat v_i\mcom\hat v_j\mcom v_{n+1})}{v_{n+2}}\\
+
\sum_{1\le i\le n+1}\s i
\ip{[v_i,\ph(v_1\mcom \hat v_i\mcom v_{n+1})]}{v_{n+2}}\\
=
\sum_{1\le i<j\le n+1}\s{i+j-1}\tph([v_i,v_j],v_1\mcom
\hat v_i\mcom\hat v_j\mcom v_{n+1},v_{n+2})\\
+
\sum_{1\le i\le n+1}\s{i+(n+2)-1}
\tph([v_i,v_{n+2}],v_1\mcom \hat v_i\mcom v_{n+1})
\\=
\sum_{1\le i<j\le n+2}\s{i+j-1}\tph([v_i,v_j],v_1\mcom
\hat v_i\mcom\hat v_j\mcom v_{n+1},v_{n+2})\\
=d\tph(v_1\mcom v_{n+2}).
\end{multline}
The last equality follows from the triviality of the action of
$\V$ on $\k$, so that the second term in the definition of the coboundary
operator drops out. Since $d\tph$ is antisymmetric, it follows that
$d\ph$ is cyclic.

A deformation $l_t=l+t\ph$ of the Lie algebra is said to
preserve the inner product if the inner product remains invariant
for $l_t$. This occurs precisely when $\ph$ is cyclic with respect to
the inner product. Similarly, a trivial deformation which preserves
the inner product is given by a linear map $\rho_t=I+t\lambda$
which satisfies
\begin{equation}
\ip{\rho_t(v_1)}{\rho_t(v_2)}=\ip{v_1}{v_2},
\end{equation}
which is equivalent to the condition
$\ip{\lambda(v_1)}{v_2}=-\ip{v_1}{\lambda(v_2)}$; in other words, $\lambda$
is cyclic.

Thus the cyclic cohomology $HC^2(\V)$ characterizes the deformations
of $\V$ which preserve the inner product. Since $HC^2(\V)\cong H^3(\V,\k)$,
we see that the cyclic cohomology is independent of the inner product.
Of course the isomorphism does depend on the inner product. Thus we
should define cyclic cohomology to be the cohomology on $C(V,\k)$
induced by the isomorphism between $CC(V)$ and $C(V,\k)$. It is this
cohomology which we have shown coincides with $H(V,\k)$, so
is independent of the inner product.

Finally, we introduce some notation that will make the definition of the
coboundary operator generalize more easily to the \zt-graded case. Let
$\sh pq$ be
the {\bf unshuffles} of type $p,q$; that is,
the subset of permutations $\sigma$ of $p+q$ elements
such that $\sigma(i)<\sigma(i+1)$ when $i\ne p$. Let $\s{\sigma}$
be the sign of the permutation. Then
\begin{multline}
d\ph(v_1\mcom v_{n+1})=
\sum_{\sigma\in\sh2{n-1}}
\s\sigma
\ph([v_{\sigma(1)},v_{\sigma(2)}],v_{\sigma(3)}\mcom v_{\sigma(n+1)})\\
-\s{n-1}
\sum_{\sigma\in\sh n1)}
\s\sigma
[\ph(v_{\sigma(1)},\mcom v_{\sigma(n)}),v_{\sigma(n+1)})]
\end{multline}

Let us introduce a \Z-grading on the \k-module $C^*(\V)=
\bigoplus_{k=1}^\infty C^k(\V)$ by
defining $\deg(\ph)=k-1$ if $\ph\in C^k(\V)$. With this grading,
$C^*(\V)$ becomes a \Z-graded Lie algebra, with the bracket defined
by
\begin{multline}
[\ph,\psi](v_1\mcom v_{k+l-1})=\\
\sum_{\sigma\in\sh l{k-1}}
\s\sigma
\ph(\psi(v_{\sigma(1)}\mcom v_{\sigma(k)}),v_{\sigma(k+1)}\mcom
v_{\sigma(k+l-1)})\\
-\s{(k-1)(l-1)}
\sum_{\sigma\in\sh k{l-1})}
\s\sigma
\psi(\ph(v_{\sigma(1)},\mcom v_{\sigma(l)}),
v_{\sigma(l+1)}\mcom v_{\sigma(n+1)}),
\end{multline}
for $\ph\in C^l(\V)$ and $\psi\in C^k(\V)$.
Later we shall establish in a more general context that this bracket
satisfies the \Z-graded Jacobi Identity. The definition of the
bracket does not depend on the fact that $\V$ is a Lie algebra,
but in the case of a Lie algebra, we easily see that
$d\ph=[\ph,l]$. That $[l,l]=0$ is an immediate consequence of the
Jacobi identity, and since $\deg(l)=1$, so that $l$ is an
odd mapping, these facts
show immediately that $d^2=0$.
A generalization of this result is that given a not
necessarily invariant inner product on $\V$,
the bracket of cyclic elements is again a cyclic element.
\begin{thm}\label{th1}
Suppose that $\V$ is a \k-module, and $\ipf$ is an
inner product on $\V$. Suppose that $\ph\in C^l(\V)$
and $\psi\in C^k(\V)$ are cyclic. Then $[\ph,\psi]$ is also
cyclic. Moreover, the formula below holds.
\begin{equation}
\widetilde{[\ph,\psi]}(v_1\mcom v_{k+l})=
\sum_{\sigma\in\sh l{k-1}}
\s{\sigma}
\tilde\ph(\psi(v_{\sigma(1)}\mcom v_{\sigma(k)}),v_{\sigma(k+1)}\mcom
v_{\sigma(k+l)}).
\end{equation}
As a consequence, the inner product induces the structure
of a graded Lie algebra in $C^*(\V,\k)$, given by
$[\tilde\ph,\tilde\psi]=\widetilde{[\ph,\psi]}$.
If further, $l:\bigwedge^2 \V\ra \V$ is a Lie algebra bracket,
and the inner product is invariant with respect to this bracket, then
the differential $d$ in $C^n(\V,\k)$ is
given by $d(\tilde\ph)=[\tilde \ph,\tilde l]$, so that
the isomorphism between $CC^*(\V)$ and $C^{*+1}(\V,\k)$, is an isomorphism
of differential $\Z$-graded Lie algebras.
\end{thm}

\section{\zt-Graded Lie Algebras}\label{sect 2}
Recall that a \zt-graded Lie algebra,
is a \zt-graded \k-module equipped with an degree zero
bracket $l:\V\bigotimes \V\ra \V$, abbreviated by $l(a,b)=[a,b]$,
which is graded anticommutative, so that
$[v_1,v_2]=\s{v_1v_2}[v_2,v_1]$, and satisfies the graded Jacobi identity
\begin{equation}
[v_1,[v_2,v_3]]=[[v_1,v_2],v_3]+\s{v_1v_2}[v_2,[v_1,v_3]].
\end{equation}
Odd brackets can also be considered,
but here we require that the bracket has degree zero, so that
$\e{[x,y]}=\e{x}+\e{y}$.
One can restrict to the case where \k\
is a field of characteristic zero, but it is also interesting to
allow \k\ to be a \zt-graded commutative ring, requiring that $\k_0$
be a field of characteristic zero.

The graded
antisymmetry of the bracket means that the bracket is a linear map
$\left[\cdot,\cdot\right]:
\bigwedge^2\V\ra \V$, where $\bigwedge^2 \V$ is the graded wedge
product. Recall that the graded exterior algebra $\bigwedge \V$ is defined
as the quotient of $\tens \V$ by the graded ideal generated by elements
of the form $x\tns y+\s{xy}y\tns x$, for homogeneous elements $x$, $y$
in $\V$. An element $v_1\wedge\dots \wedge v_n$ is said to have
(external) degree
$n$ and parity (internal degree)
$\e{v_1}\mplus \e{v_n}$. If $\omega$ has degree
$\deg(\omega)$ and parity $\e{\omega}$, and similarly for $\eta$, then
\begin{equation}
\omega \wedge \eta=\s{\e\omega\e\eta+\deg(\omega)\deg(\eta)}
\eta\wedge\omega.
\end{equation}
More generally (see section \ref{sect 3}),
if $\sigma$ is any permutation, we define
$\epsilon(\sigma;v_1\mcom v_n)$ by requiring that
\begin{equation}
v_1\mwedge v_n=\s{\sigma}
\epsilon(\sigma;v_1\mcom v_n)
v_{\sigma(1)}\mwedge v_{\sigma(n)},
\end{equation}
where $\s{\sigma}$ is the sign of the permutation $\sigma$.

In order to see how the coboundary operator should be modified in the case
of \zt-graded algebras, we consider infinitesimal deformations of
the Lie algebra. If we denote the deformed bracket by $l_t=l+t\ph$ as
before, then we wish the deformed bracket to remain even, so
that if $t$ is taken to be an even parameter,
then $\ph$ must be even.
However, if we let $t$ be an odd parameter, then $\ph$ must
be odd. We also must take into account that parameters should
graded commute with elements of \V,
so that $vt=\s{ta}tv$.
The graded Jacobi identity takes the form
\begin{equation}
l_t(a,l_t(b,c))=l_t(l_t(a,b),c)+\s{ab}l_t(b,l_t(a,c))
\end{equation}
Expanding this formula, we obtain the condition
\begin{multline}
\ph(l(a,b),c)+\s{bc+1}\ph(l(a,c),b)+\s{a(b+c)}\ph(l(b,c),a)+\\
l(\ph(a,b),c)+\s{bc+1}l(\ph(a,c),b)+\s{a(b+c)}l(\ph(b,c),a)=0,
\end{multline}
which can be expressed as $d\ph=0$, if we define for
$\ph:\bigwedge^n \V\ra \V$,
\begin{multline}\label{ztcb}
d\ph(v_1\mcom v_{n+1})=
\sum_{\sigma\in\sh2{n-1}}
\s\sigma\epsilon(\sigma)
\ph([v_{\sigma(1)},v_{\sigma(2)}],v_{\sigma(3)}\mcom v_{\sigma(n+1)})\\
-\s{n-1}
\sum_{\sigma\in\sh n1)}
\s\sigma\epsilon(\sigma)
[\ph(v_{\sigma(1)},\mcom v_{\sigma(n)}),v_{\sigma(n+1)})].
\end{multline}
More generally, when $M$ is a
\zt-graded $\V$ module, we can define a right
multiplication $M\tens \V\ra M$ by $m\cdot a=-\s{ma}a\cdot m$. Then
if we let $C^n(\V,M)=\hom(\bigwedge^n\V,M)$ be the module of
cochains of degree $n$ on \V\ with values in $M$, we can define
a coboundary operator $d:C^n(\V,M)\ra C^{n+1}(\V,M)$ by
\begin{multline}
d\ph(v_1\mcom v_{n+1})=
\sum_{\sigma\in\sh2{n-1}}
\s\sigma\epsilon(\sigma)
\ph([v_{\sigma(1)},v_{\sigma(2)}],v_{\sigma(3)}\mcom v_{\sigma(n+1)})\\
-\s{n-1}
\sum_{\sigma\in\sh n1)}
\s\sigma
\ph(v_{\sigma(1)},\mcom v_{\sigma(n)})\cdot v_{\sigma(n+1)}.
\end{multline}
Then we define $H^n(\V,M)$ in the same manner as for
ordinary Lie algebras, and as before, denote $C^n(\V,\V)=C^n(\V)$
and $H^n(\V)=H(V,V)$ for the adjoint action of \V.

As in the case of ordinary Lie algebras, we obtain that trivial
(infinitesimal) deformations
are given by means of coboundaries of linear maps $\lambda:\V\ra \V$.
Thus the infinitesimal deformations of a \zt-graded Lie algebra
are classified by $H^2(\V)$, where the cohomology of the Lie algebra is
given by means of the coboundary operator in equation (\ref{ztcb})
above. The fact that $d^2=0$ is verified in the usual manner.

$C^*(\V)=\bigoplus_{k=1}^\infty C^*(\V)$ has a natural $\zt\times\Z$
grading, with the bidegree of a homogeneous element $\ph\in C^n(\V)$
given by $\bid(\ph)=(\e\ph,n-1)$. (Actually, it is more common to
define the exterior degree of $\ph$ to be $1-n$, but as the degree
enters into our calculations mainly by sign, this makes no difference.)
There is a natural bracket operation
in $C^*(\V)$, given by
\begin{multline}\label{nbra}
[\ph,\psi](v_1\mcom v_{k+l-1})=\\
\sum_{\sigma\in\sh l{k-1}}\!\!
\s\sigma\epsilon(\sigma)
\ph(\psi(v_{\sigma(1)}\mcom v_{\sigma(k)}),v_{\sigma(k+1)}\mcom
v_{\sigma(k+l-1)})\\
-\s{\ph\psi+ (k-1)(l-1)}\!\!\!\!\!\!\!\!\!\!
\sum_{\sigma\in\sh k{l-1})}\!\!
\s\sigma\epsilon(\sigma)
\psi(\ph(v_{\sigma(1)},\mcom v_{\sigma(l)}),
v_{\sigma(l+1)}\mcom v_{\sigma(n+1)}).
\end{multline}
This bracket makes $C^*(\V)$ into a $\zt\times\Z$-graded Lie algebra.
The differential takes the form $d(\ph)=[\ph,l]$, and the condition
$[l,l]=0$ is precisely equivalent to the \zt-graded Jacoby identity
for $l$. Since $l$ has even parity, we note that as a
$\zt\times\Z$-graded map, it is odd, since $\ip{(0,1)}{(0,1)}=1$.

For a \zt-graded \k-module \V, an inner product is
a (right) \k-module
homomorphism $h:\V\tens\V\ra\V$, which is (graded) symmetric, and
non-degenerate. Denote the inner product by $\ip vw=h(v\tns w)$.
Graded symmetry means that $\ip vw=\s{vw}\ip wv$.
Non-degeneracy means that the map $\lambda:\V\ra \V^*=\hom(\V,\k)$, given by
$\lambda(v)(w)=\ip vw$ is an isomorphism. (When \k\ is a field, this is
equivalent to the usual definition of a non degenerate bilinear form.)
If $h$ is an even map, then we say that the inner product is even. We
shall only consider even inner products on \V.
If $\V$ is a \zt-graded Lie algebra,
then we define the notion of an invariant inner product
as in the non graded case, by
$\ip{[v_1,v_2]}{v_3}=\ip{v_1}{[v_2,v_3]}$.
Then the tensor $\tilde l$, given by
\begin{equation}
\tilde l(v_1,v_2,v_3)=\ip{\sb{v_1,v_2}}{v_3},
\end{equation}
is (graded) antisymmetric, so that $\tilde l\in\hom(\bigwedge^3\V,\k)$.
In particular, we have
\begin{equation}
\tilde l(v_1,v_2,v_3)=
\s{v_3(v_1+v_2)}
\tilde l(v_3,v_1,v_2),
\end{equation}
so that $\tilde l$ satisfies a property of graded invariance under
cyclic permutations. In general, we say that an element $\ph\in
C^n(\V)=\hom(\bigwedge^n\V,\V)$ is cyclic with respect to the inner
product, or preserves the inner product, if
\begin{equation}
\ip{\ph(v_1\mcom v_n)}{v_{n+1}}=\s{n+v_1\ph}
\ip{v_1}{\ph(v_2\mcom v_{n+1})}.
\end{equation}
Then $\ph$ is cyclic if and only if $\tilde\ph:\bigwedge^{n+1}\V
\ra \V$, given by
\begin{equation}
\tilde\ph(v_1\mcom v_{n+1})=
\ip{\ph(v_1\mcom v_n)}{v_{n+1}}
\end{equation}
is antisymmetric. Since the inner product is non-degenerate,
the map $\ph\mapsto\tilde\ph$ is an even isomorphism between
the submodule $CC^n(\V)$ of $C^n(\V)$ consisting of cyclic
elements, and $C^{n+1}(\V,\k)$.

We obtain a straightforward generalization of theorem
\ref{th1} in the \zt-graded case.
\begin{thm}\label{th2}
Suppose that $\V$ is a \zt-graded \k-module, and $\ipf$ is an
inner product on $\V$. Suppose that $\ph\in C^l(\V)$
and $\psi\in C^k(\V)$ are cyclic. Then $[\ph,\psi]$ is also
cyclic. Moreover, the formula below holds.
\begin{equation}
\widetilde{[\ph,\psi]}(v_1\mcom v_{k+l})=\!\!\!\!\!
\sum_{\sigma\in\sh l{k-1}}\!\!\!\!\!
\s{\sigma}\epsilon(\sigma)
\tilde\ph(\psi(v_{\sigma(1)}\mcom v_{\sigma(k)}),v_{\sigma(k+1)}\mcom
v_{\sigma(k+l)}).
\end{equation}
As a consequence, the inner product induces the structure
of a \ztz-graded Lie algebra in $C^*(\V,\k)$, given by
$[\tilde\ph,\tilde\psi]=\widetilde{[\ph,\psi]}$.
If further, $l:\bigwedge^2 \V\ra \V$ is a
\zt-graded Lie algebra bracket,
and the inner product is invariant with respect to this bracket, then
the differential $d$ in $C^n(\V,\k)$ is
given by $d(\tilde\ph)=[\tilde \ph,\tilde l]$, so that
the isomorphism between $CC^*(\V)$ and $C^{*+1}(\V,\k)$, is an isomorphism
of differential $\zt\times\Z$-graded Lie algebras.
\end{thm}
In particular, we can define cyclic cohomology of a \zt-graded Lie
algebra in the same manner as for an ordinary Lie algebra. The
isomorphism $HC^n(\V)\cong H^{n+1}(\V,\k)$ holds for the \zt-graded case
as well.

In the following sections, we shall generalize the notion of cohomology
of a graded Lie algebra, and cyclic cohomology of the same, to the
case of \linf\ algebras. The important observation is that cohomology
of a Lie algebra is determined by a bracket operation on the cochains.
We shall see that this leads to the result that cohomology of a Lie
algebra classifies the deformations of the Lie algebra into
an \linf\ algebra.
\section{The Exterior and Symmetric Algebras}\label{sect 3}
Suppose that $V$ is a \zt-graded
\k-module. For the purposes of this paper, we shall define the tensor
algebra $T(V)$ by
$T(V)=\bigoplus_{n=1}^\infty V^n$, where $V^n$ is the $n$-th tensor
power of $V$. Often, the tensor algebra is defined to include the
term $V^0=\k$ as well, but we shall omit it here. For a treatment of
\ainf\ algebras which includes this term see \cite{getz}.
For an element $v=v_1\mtns v_n$ in $T(V)$, define
its parity $\e{v}=\e{v_1}\mplus \e{v_n}$, and its degree by
$\deg(v)$=n. We define the bidegree of $v$ by
$\bd v=(\e{v},\deg(v))$. If $u, v \in T(V)$, then $\bd{u\tns v}=
\bd u +\bd v$, so that $T(V)$ is naturally $\zt\times\Z$-graded by
the bidegree, and \zt-graded if we consider only the parity.
If $\sigma\in\Sigma_n$, then there is a natural right \k-module
homomorphism $S_\sigma:V^n\ra V^n$, which satisfies
\begin{equation}\label{a1}
S_\sigma(v_1\mtns v_n)=\epsilon(\sigma;v_1\mcom v_n)
v_{\sigma(1)}\mtns v_{\sigma(n)},
\end{equation}
where $\epsilon(\sigma;v_1\mcom v_n)$ is a sign which can be determined
by the following. If $\sigma$ interchanges $k$ and $k+1$, then
$\epsilon(\sigma;v_1\mcom v_n)=\s{v_kv_{k+1}}$. In addition,
if $\tau$ is another permutation, then
\begin{equation}
\epsilon(\tau\sigma;v_1\mcom v_n)=
\epsilon(\tau;v_{\sigma(1)}\mcom v_{\sigma(n)})
\epsilon(\sigma;v_1\mcom v_n).
\end{equation}
It is conventional to abbreviate $\epsilon(\sigma;v_1\mcom v_n)$
as $\epsilon(\sigma)$.

The symmetric algebra $\bigodot V$ is defined as the quotient
of the tensor algebra of $V$ by the bigraded ideal generated by all
elements of the form $u\tns v-\s{uv}v\tns u$. The resulting algebra
has a decomposition $\bigodot V=\bigoplus_{n=1}^\infty \bigodot^n V$,
and the induced product is denoted by $\odot$.
The symmetric algebra is both \zt\ and $\zt\times\Z$-graded, and is graded
commutative with respect to the \zt-grading. For simplicity, let us
denote $\s{uv}=\s{\e u\e v}$.
In other words,
if $u, v\in\bigodot V$, then
\begin{equation}
u\odot v=\s{uv}v\odot u.
\end{equation}
Furthermore, it is easy to see that
\begin{equation}
v_1\modot v_n=\epsilon(\sigma)v_{\sigma(1)}\modot v_{\sigma(n)}.
\end{equation}

The exterior algebra $\bigwedge V$ is defined
as the quotient of
$T(V)$ by the bigraded ideal generated by all
elements of the form $u\tns v+\s{uv}v\tns u$. The resulting algebra
has a decomposition $\bigwedge V=\bigoplus_{n=0}^\infty \bigwedge^n V$,
and the induced product is denoted by $\wedge$.
The exterior algebra is both \zt\ and $\zt\times\Z$-graded, and is
graded commutative with respect to the $\zt\times\Z$-grading.
We introduce a \zt-valued
inner product on $\zt\times\Z$ by
\begin{equation}
\ip{(\bar m,n)}{(\bar r,s)}=
{\bar m\bar r+\bar n\bar s}.
\end{equation}
Let $u, v\in \bigwedge V$. For simplicity, let us denote
$\s{\ip uv}=\s{\ip{\bid(u)}{\bid(v)}}$.
Then
$$u\wedge v=\s{\ip uv} v\wedge u,$$
which is precisely the formula for $\zt\times\Z$-graded
commutativity. Furthermore it is easy to see that
\begin{equation}
v_1\mwedge v_n=\s{\sigma}\epsilon(\sigma)v_{\sigma(1)}\mwedge v_{\sigma(n)}.
\end{equation}

\section{The Tensor, Exterior and Symmetric Coalgebras}\label{sect
4}
A proper treatment of the symmetric and exterior coalgebras would
introduce the coalgebra structure on the tensor algebra, and then
describe these coalgebras in terms of coideals. Instead, we will
describe these coalgebra structures directly. Recall that a coalgebra
structure on a \k-module $C$ is given by a diagonal mapping
$\Delta:C\ra C\tens C$. We consider only coassociative coalgebras,
but we do not consider counits. The axiom of coassociativity is
that $(1\tns \Delta)\circ \Delta=(\Delta\tns 1)\circ\Delta$.
A grading on $C$ is compatible  with the coalgebra structure
if for homogeneous $c\in C$, $\Delta(c)=\sum_i u_i\tens v_i$, where
$\e{u_i}+\e{v_i}=\e c$ for all $i$.
We also mention that a coalgebra is graded cocommutative if
$S\circ\Delta=\Delta$, where $S:C\tens C\ra C\tens C$ is the
symmetric mapping given by $S(m\tns n)=\s{\ip mn}n\tns m$.

The tensor coalgebra structure is given by defining the (reduced) diagonal
$\Delta:T(V)\ra T(V)$ by
\begin{equation}
\Delta(v_1\mtns v_n)=\sum_{k=1}^{n-1}(v_1\mtns v_k)\tns(v_{k+1}\mtns v_n).
\end{equation}
(We use here the reduced diagonal, because we are not including the zero
degree term in the tensor coalgebra.)
The tensor coalgebra is not graded cocommutative under either the \zt\
or the $\zt\times\Z$-grading, but both gradings are compatible with
the coalgebra structure.

The symmetric coalgebra structure on $\bigodot V$ is given by
defining
\begin{equation}
\Delta(v_1\modot v_n)=
\sum_{k=1}^{n-1}
\sum_{\sigma\in\sh k{n-k}}
\epsilon(\sigma)
v_{\sigma(1)}\modot v_{\sigma(k)}
\tns
v_{\sigma(k+1)}\modot v_{\sigma(n)}.
\end{equation}
With this coalgebra structure, and the \zt-grading, $\bigodot V$
is a cocommutative, coassociative coalgebra without a counit.
Similarly, we define the exterior coalgebra structure on
$\bigwedge V$ by
\begin{equation}
\Delta(v_1\mwedge v_n)=
\sum_{k=1}^{n-1}
\sum_{\sigma\in\sh k{n-k}}\!\!\!\!\!\!
\s\sigma
\epsilon(\sigma)
v_{\sigma(1)}\mwedge v_{\sigma(k)}
\tns
v_{\sigma(k+1)}\mwedge v_{\sigma(n)}.
\end{equation}
Then the coalgebra structure is coassociative, and is cocommutative
with respect to the $\zt\times\Z$-grading.
\section{Coderivations}\label{sect 5}
A coderivation on a graded coalgebra $C$ is a map $d:C\ra C$
such that
\begin{equation}
\Delta\circ d=(d\tns 1+1\tns d)\circ\Delta.
\end{equation}
Note that the definition depends on the grading group, because
$(1\tns d)(\alpha\tns\beta)=\s{\alpha d}\alpha\tns d(\beta)$.
The \k-module $\coder(C)$
of all graded coderivations has a natural structure
of a graded Lie algebra, with the bracket given by
\begin{equation}
\sb{m,n}=m\circ n-\s{\ip mn}n\circ m,
\end{equation}
where $\s{\ip mn}=\s{\ip{\e m}{\e n}}$ is given by the inner product
in the grading group, so that the definition of the bracket also
depends on the grading group.
A codifferential on a coalgebra $C$ is a coderivation $d$ such that
$d\circ d=0$.
We examine the
coderivation
structure of the tensor, symmetric and exterior coalgebras.
\subsection{Coderivations of the Tensor Coalgebra}
Suppose that we wish to extend $d_k:V^k\ra V$ to a coderivation
$\hd_k$ of $T(V)$.
We are interested in extensions satisfying the property that
$\hd_k(v_1\mcom v_n)=0$ for $n<k$.
How this extension is made depends on whether we consider the \zt\ or
the $\zt\times\Z$-grading. First we consider the \zt-grading, so that
only the parity of $d$ is relevant. Then if we define
\begin{multline}
\hd_k(v_1\mtns v_n)=\\
\sum_{i=0}^{n-k}
\s{(v_1\mplus v_i)d_k}
v_1\mtns v_i\tns d_k(v_{i+1}\mcom v_{i+k})\tns v_{i+k+1}\mtns v_n,
\end{multline}
one can show that $\hd_k$ is a coderivation on $T(V)$ with
respect to the \zt-grading.
More generally, one can show that any coderivation $\hd$ on $T(V)$ is
completely determined by the induced mappings $d_k:V^k\ra V$, and in
fact, one obtains that
\begin{multline}\label{cder1}
\hd(v_1\mtns v_n)=\\
\sum
\begin{Sb}
1\le k\le n\\
\\
0\le i\le n-k
\end{Sb}
\s{(v_1\mplus v_i)d_k}
v_1\mtns v_i\tns d_k(v_{i+1}\mcom v_{i+k})\tns v_{i+k+1}\mtns v_n.
\end{multline}
Also, one can show that $\hd$ is a codifferential with respect to the
\zt-grading precisely when
\begin{equation}
\sum
\begin{Sb}
k+l=n+1\\
\\
0\le i\le n-k
\end{Sb}
\s{(v_1\mplus v_i)d_k}
d_l(v_1\mcom v_i,d_k(v_{i+1}\mcom v_{i+k}),v_{i+k+1}\mcom v_n)=0.
\end{equation}
The module
$\coder(T(V))$ of coderivations of $T(V)$ with respect to the
\zt-grading
is naturally isomorphic to $\hom(T(V),V)$, so $\hom(T(V),V)$ inherits
a natural structure of a \zt-graded Lie algebra.

Let us examine the bracket structure on $\hom(T(V),V)$ more closely.
Suppose that for an arbitrary element $d\in\hom(T(V),V)$,
we denote by $d_k$ the restriction of $d$ to $V^k$, and by
$\hd$, $\hd_k$  the extensions of $d$ and $d_k$
as coderivations of $T(V)$. Also denote by $d_{kl}$ the
restriction of of $\hd_k$ to $V^{k+l-1}$,
so that $d_{kl}\in\hom(V^{k+l-1},V^l)$.
The precise expression for $d_{kl}$ is given by equation (\ref{cder1})
with $n=k+l-1$. It is easy to see that the bracket of $d_k$ and
$\delta_l$ is given by
\begin{equation}
[d_k,\delta_l]=
d_k\circ\delta_{lk}-\s{d_k\delta_l}\delta_l\circ d_{kl}.
\end{equation}
Furthermore, we have
$[d,\delta]_n=\sum_{k+l=n+1}[d_k,\delta_l].$
The point here is that $d_{kl}$ and $\delta_{kl}$ are determined in a
simple manner by $d_k$ and $\delta_k$, so we have given a description
of the bracket on $\hom(T(V),V)$ in a direct fashion. The fact that
the bracket so defined has the appropriate properties follows from
the fact that if $\rho=[d,\delta]$, then $\hat\rho=[\hd,\hdl]$.

Now we consider how to extend a mapping $m_k:V^k\ra V$ to a coderivation
$\hm_k$ with respect to the \ztz-grading. In this case, note that the
bidegree of $m_k$ is given by $\bd{m_k}=(\e{m_k},k-1)$. The formula for the
extension is the same as before, but with the bidegree in place of
the parity. In other words,
\begin{multline}
\hm_k(v_1\mtns v_n)=\\
\sum_{i=0}^{n-k}
\s{(v_1\mplus v_i)m_k+i(k-1)}
v_1\mtns v_i\tns m_k(v_{i+1}\mcom v_{i+k})\tns v_{i+k+1}\mtns v_n.
\end{multline}
Similarly, if we consider an arbitrary coderivation $\hm$ on $T(V)$,
then it is again determined by the induced mappings $m_k:V^k\ra V$,
and we see that
\begin{multline}
\hm(v_1\mtns v_n)=
\sum
\begin{Sb}
1\le k\le n\\
\\
0\le i\le n-k
\end{Sb}
\s{(v_1\mplus v_i)m_k +i(k-1)}\\\times
v_1\mtns v_i\tns m_k(v_{i+1}\mcom v_{i+k})\tns v_{i+k+1}\mtns v_n.
\end{multline}
Also, one obtains that $\hm$ is a codifferential with respect to the
$\zt\times\Z$-grading is equivalent to the condition that for all $n$,
\begin{multline}\label{ztzcodiff}
\sum
\begin{Sb}
k+l=n+1\\
\\
0\le i\le n-k
\end{Sb}
\s{(v_1\mplus v_i)m_k +i(k-1)}\\\times
m_l(v_1\mcom v_i,m_k(v_{i+1}\mcom v_{i+k}),v_{i+k+1}\mcom v_n)=0.
\end{multline}

The module
$\coder(T(V))$ of coderivations of $T(V)$ with respect to the
\ztz-grading
is naturally isomorphic to
$\bigoplus_{k=1}^\infty\hom(V^k,V)$, rather than $\hom(T(V),V)$, because
the latter module is the direct product of the modules $\hom(V^k,V)$.
However, we would like to consider
elements of the form $\hm=\sum_{k=1}^\infty\hm_k$, where $\hm_k$
has bidegree $(\e{m_k},k-1)$. Such an infinite sum is a well
defined element of $\hom(T(V),T(V))$, so by abuse of notation, we
will define $\coder(T(V))$ to be the module of such infinite sums
of coderivations. With this convention, we now have a natural
isomorphism between $\coder(T(V))$ and $\hom(T(V),V)$. Furthermore,
the bracket of coderivations is still well defined, and we consider
$\coder(T(V)$ to be a \ztz-graded Lie algebra. The reason that the
bracket is well defined is that any homogeneous coderivation has
bidegree $(m,n)$ for some $n\ge 0$, so the grading is
given by $\zt\times\N$ rather than the full group \ztz. In structures
where a $\Z$-grading reduces to an \N-grading, it is often advantageous
to replace direct sums with direct products.

Using the same notation convention as in the \zt-graded case, we
note that if $m,\mu\in\hom(T(V),V)$, then we have
\begin{equation}
[m_k,\mu_l]=m_k\circ\mu_{lk}-\s{\ip{m_k}{\mu_l}}\mu_l\circ m_{kl},
\end{equation}
and $[m,\mu]_n=\sum_{k+l=n+1}[m_k,\mu_l]$.

\subsection{Coderivations of the Symmetric Coalgebra}
Suppose that we want to extend $m_k:\bigodot^k V\ra V$ to a coderivation
$\hm_k$ of $\bigodot V$ such that
$m(v_1\modot v_n)=0$ for $k<n$.
Define
\begin{equation}
\hm_k(v_1\modot v_n)=
\sum_{\sigma\in\sh k{n-k}}
\epsilon(\sigma)
m_k(v_{\sigma(1)}\mcom v_{\sigma(k)})\odot v_{\sigma(k+1)}
\modot v_{\sigma(n)}.
\end{equation}
Then $\hm_k$ is a coderivation with respect to the \zt-grading.
In general, suppose that $\hm$ is a coderivation on the  symmetric
coalgebra. It is not difficult to see that if
$m_k:\bigodot^k V\ra V$ is the induced map,
then $\hm$ can be recovered
from these maps by the relations
\begin{equation}
\hm(v_1\modot v_n)=
\sum
\begin{Sb}
1\le k\le n\\
\\
\sigma\in\sh k{n-k}
\end{Sb}
\epsilon(\sigma)
m_k(v_{\sigma(1)}\mcom v_{\sigma(k)})\odot v_{\sigma(k+1)}
\modot v_{\sigma(n)}.
\end{equation}
{}From this, we determine that there is a natural isomorphism
between
$\coder(\bigodot V)$, the module of coderivations of $V$,
and $\hom(\bigodot V,V)$.
Thus $\hom(\bigodot V,V)$ inherits
the structure of a graded Lie algebra.
Also, $\hm$ is a codifferential when for all $n$,
\begin{equation}
\sum
\begin{Sb}
k+l=n+1
\\
\sigma\in\sh k{n-k}
\end{Sb}
\epsilon(\sigma)
m_l(m_k(v_{\sigma(1)}\mcom v_{\sigma(k)}),v_{\sigma(k+1)}
\mcom v_{\sigma(n)})=0.
\end{equation}

It is reasonable to ask whether a map
$m_k:\bigodot^k V\ra V$ can be extended as a coderivation with respect
to the \ztz-grading. It turns out that in general, it is not possible
to do this. For example, suppose that we are given
$m_2:\bigodot^2 V\ra V$. If $m_2$ is extendible as a coderivation $\hm$,
then we must have
\begin{multline}
\Delta \hm_2(v_1\odot v_2\odot v_3)=(\hm\tns 1+1\tns \hm)
\Delta
(v_1\odot v_2\odot v_3)=\\
=(\hm\tns 1+1\tns \hm)[v_1\tns v_2\odot v_3 +\s{v_1v_2}v_2\tns v_1\odot
v_3+\\
\s{v_3(v_1+v_2)}v_3\tns v_1\odot v_2+ v_1\odot v_2\tns v_3+\\
\s{v_2v_3}v_1\odot v_3\tns v_2+\s{v_1(v_2+v_3)}v_2\odot v_3\tns v_1]\\
=\s{\ip{v_1}{m_2}}v_1\tns m_2(v_2,v_3)+\s{\ip{v_2}{m_2}+v_1v_2}
v_2\tns m_2(v_1,v_3)+\\
\s{\ip{v_3}{m_2}+v_3(v_1+v_2)}v_3\tns m_2(v_1,v_2)+ m_2(v_1,v_2)\tns v_3 +\\
\s{v_2v_3}m_2(v_1,v_3)\tns v_3+\s{v_1(v_2+v_3)}m_2(v_2,v_3)\tns v_1
\end{multline}
In the above we are using the fact that
$(1\tns m_2)(\alpha\tns \beta)=\s{\ip \alpha{m_2}}a\tns m_2(\beta)$,
where $\ip\alpha {m_2}$ depends on which grading group we are using.
We need to recognize the expression as $\Delta$ of something. In order
for this to be possible, the terms $\s{\ip{v_1}{m_2}}v_1\tns m_2(v_2,v_3)$
and $\s{v_1(v_2+v_3)}m_2(v_2,v_3)\tns v_1$ need to match up. We have
$\Delta(v_1\tns m_2(v_2,v_3))=\s{v_1(v_2+v+3+m_2)}m_2(v_2,v_3)\tns v_1$.
Thus for the expressions to match up, it is necessary that
$\ip{v_1}{m_2}=\e{v_1}\e{m_2}$, which is the inner product given by the
\zt-grading.
\subsection{Coderivations of the Exterior Coalgebra}
Suppose that we want to extend $l_k:\bigwedge^k V\ra V$ to a
coderivation $\hl_k$ of $\bigwedge V$.
Define
\begin{equation}
\hl(v_1\mwedge v_n)=
\sum_{\sigma\in\sh k{n-k}}
\s{\sigma}
\epsilon(\sigma)
l_k(v_{\sigma(1)}\mcom v_{\sigma(k)})\wedge v_{\sigma(k+1)}
\mwedge v_{\sigma(n)}.
\end{equation}
Then $\hl$ is a coderivation with respect to the $\zt\times\Z$-grading.
As in the graded symmetric case, an arbitrary coderivation on
$\bigwedge V$ is completely determined by the induced maps
$l_k:\bigwedge^k V\ra V$ by the formula
\begin{equation}
l(v_1\mwedge v_n)=
\sum
\begin{Sb}
1\le k\le n\\
\\
\sigma\in\sh k{n-k}
\end{Sb}
\s{\sigma}
\epsilon(\sigma)
l_k(v_{\sigma(1)}\mcom v_{\sigma(k)})\wedge v_{\sigma(k+1)}
\mwedge v_{\sigma(n)}.
\end{equation}
Similarly to the previous cases considered, we see that
$\hl$ is a codifferential when for all $n$,
\begin{equation}
\sum
\begin{Sb}
k+l=n+1\\
\\
\sigma\in\sh k{n-k}
\end{Sb}
\s{\sigma}
\epsilon(\sigma)
l_l(l_k(v_{\sigma(1)}\mcom v_{\sigma(k)})\wedge v_{\sigma(k+1)}
\mwedge v_{\sigma(n)})=0.
\end{equation}
As in the case of the tensor coalgebra, the \ztz-grading requires
us to extend the notion of a coderivation in order to obtain an
isomorphism between $\coder(\bigwedge V)$ and $\hom(\bigwedge V,V)$.
In our extended sense, $\hom(\bigwedge V,V)$ inherits a natural
structure of a \ztz-graded Lie algebra.

One can ask whether a map $l_k:\bigwedge^k V\ra V$ can be extended
as a coderivation with respect to the \zt-grading.
For example, suppose that $l:\bigwedge^2 V\ra V$ is given. Then
if $l$ is to be extended as a coderivation with respect to the
$\zt\times\Z$-grading, we must have
\begin{multline}
\Delta l(v_1\wedge v_3\wedge v_3)=
(l\tns1+1\tns l)
\!\!\!\!
\sum
\begin{Sb}
1\le k\le 2\\
\\
\sigma\in\sh k{3-k}
\end{Sb}
\!\!\!\!
\s{\sigma}
\epsilon(\sigma)
v_{\sigma(1)}\mwedge v_{\sigma(k)}
\tns
v_{\sigma(k+1)}\mwedge v_{\sigma(3)}
\\
=(l\tns1+1\tns l)[
v_1\tns v_2\wedge v_3+v_1\wedge v_2\tns v_3
+\s{v_1v_2+1}v_2\tns v_1\wedge v_3 +\\
\s{v_2v_3+1}v_1\wedge v_3\tns v_2+
\s{v_1(v_2+v_3)}v_2\wedge v_3\tns v_1+
\s{v_3(v_1+v_2)}v_3\tns v_1\wedge v_2]\\
=\s{v_1l+1}v_1\tns l(v_2,v_3)+
l(v_1,v_2)\tns v_3+
\s{v_1v_2+v_2l}v_2\tns l(v_1,v_3)+\\
\s{v_2v_3+1}l(v_1,v_3)\tns v_2+
\s{v_1(v_2+v_3)}l(v_2,v_3)\tns v_1+
\s{v_3(v_1+v_3)+v_3l+1}v_3\tns l(v_1,v_2)\\
=\Delta[l(v_1,v_2)\wedge v_3+\s{v_2v_3+1}l(v_1,v_3)\wedge v_2
+\s{v_1(v_2+v_3)}l(v_2,v_3)\wedge v_1]
\end{multline}
The map $l$ above cannot be extended as a coderivation with respect to
the \zt-grading. The signs introduced by the exchange rule when
applying $(1\tns l)$ would make it impossible to express the result
as $\Delta$ of something. Thus we need the $\zt\times\Z$-grading to
obtain a good theory of coderivations of the exterior coalgebra. Similarly,
the \zt-grading is necessary to have a good theory of coderivations of
the symmetric coalgebra.
\section{Cohomology of \ainf\ algebras}\label{sect 6}
In \cite{ps2}, a generalization of an associative algebra, called
a strongly homotopy associative algebra, or \ainf\ algebra was
discussed, and cohomology and cyclic cohomology of this structure
was defined. \ainf\ algebras were introduced by
J. Stasheff in \cite{sta1,sta2}. An \ainf\ algebra
structure is simply a codifferential on the tensor coalgebra; an
associative algebra is simply a codifferential determined by a
single map $m_2:V^2\ra V$. We present a description of the basic
theory of \ainf\ algebras, from the coalgebra point of view. Hopefully,
this will make the presentation of \linf\ algebras, which are given
by codifferentials on the exterior coalgebra, seem more natural.

If $V$ is a \zt-graded \k-module, then the parity reversion
$\Pi V$ is the same module, but with the parity of elements reversed.
In other words, $(\Pi V)_0=V_1$ and $(\Pi V)_1=V_0$, where $V_0$ and
$V_1$ are the submodules of even and odd elements of $V$, resp. The
map $\pi:V\ra \Pi V$, which is the identity as a map of sets, is
odd.
There is a natural isomorphism $\eta:T(V)\ra T(\Pi V)$ given by
\begin{equation}
\eta(v_1\mtns v_n)=
\s{(n-1)v_1\mplus v_{n-1}}\pi v_1\mtns \pi v_n
\end{equation}
Denote the restriction of $\eta$ to $V^k$ by $\eta_k$. Note that
$\eta_k$ is odd when $k$ is odd and even when $k$ is even, so that
$\eta$ is neither an odd nor an even mapping.

Let $W=\Pi V$. Define a bijection between $C(W)=\hom(T(W),W)$ and
$C(V)=\hom(T(V),V)$ by setting $\mu=\eta^{-1}\circ\delta\circ\eta$,
for $\delta\in\CW$.
Then $\mu_k=\eta_1^{-1}\circ\delta_k\circ\eta_k$ and
$\e{\mu_k}=\e{\delta_k}+(k-1)$. In particular, note that if $\delta_k$
is odd in the \zt-grading, then $\bid(\mu_k)=(k,k-1)$, so that
$\mu_k$ is odd in the \ztz-grading.
Now
extend $\delta_k:W^k\ra W$ to a coderivation $\hdl_k$
on $T(W)$ with respect to the \zt\ grading, so that
\begin{multline}
\hdl_k(w_1\mtns w_n)=
\sum_{i=0}^{n-k}
\s{(w_1\mplus w_i)\delta_k}\\\times
w_1\mtns w_i
\tns
\delta_k(w_{i+1}\mcom w_{i+k})
\tns
w_{i+k+1}\mtns w_n.
\end{multline}
Let $\bmu_k:T(V)\ra T(V)$ be given by
$\bmu_k=\eta^{-1}\circ \hdl_k\circ \eta$.
Let $\hmu$ be the extension of $\mu$ as a \ztz-graded coderivation
of $T(V)$.
We wish to investigate the
relationship between $\bmu_k$ and $\hmu_k$.
For simplicity, write $w_i=\pi v_i$.
Note that $\eta_1=\pi$ is the parity reversion
operator.  So we have $\delta_k=\pi\circ\mu_k\circ\eta_k\inv$.
Thus
\begin{multline}
\bmu_k(v_1\mcom v_n)=
\s{r}
\eta^{-1}\delta(w_1\mcom w_n)=\\
\eta^{-1}(
\sum_{i=0}^{n-k}
\s{r+s_i}
w_1\mtns w_i\tns \delta_k(w_{i+1}\mcom w_{i+k})\tns w_{i+k+1}\mtns w_{n}
)
=\\
\eta^{-1}(
\sum_{i=0}^{n-k}
\s{r+s_i +t_i}
w_1\mtns w_i\tns \pi \mu_k(v_{i+1}\mcom v_{i+k})\tns w_{i+k+1}\mtns w_{n}
)
=\\
\sum_{i=0}^{n-k}
\s{r+s_i +t_i+u_i}
v_1\mtns v_i\tns \mu_k(v_{i+1}\mcom v_{i+k})\tns v_{i+k+1}\mtns v_{n},
\end{multline}
where
\begin{eqnarray*}
r&=&(n-1)v_1\mplus v_{n-1}\\
s_i&=&(w_1\mplus w_i)\delta_k\\
&=&(v_1\mplus v_i)\mu_k+i(\mu_k+1-k)
+(1-k)(v_1\mplus v_i)\\
t_i&=&(k-1)v_{i+1}\mplus v_{i+k-1}\\
u_i&=&(n-k)v_1\mplus (n-k-i+1)v_i\\
&&+(n-k-i)(\mu_k+v_{i+1}\mplus v_{i+k})
+ (n-k-i-1)v_{i+k+1}\mplus v_{n-1}
\end{eqnarray*}
Combining these coefficients we find that
\begin{equation}
r+s_i+t_i+u_i=(v_1\mplus v_i)\mu_k+i(k-1) +(n-k)\mu_k,
\end{equation}
so that
\begin{multline}\label{hass2}
\bmu_k(v_1\mcom v_n)=
\sum_{i=0}^{n-k}
\s{ (v_1\mplus v_i)\mu_k+i(k-1) +(n-k)\mu_k}\\\times
v_1\mtns v_i\tns \mu_k(v_{i+1}\mcom v_{i+k})\tns v_{i+k+1}\mtns v_{n}.
\end{multline}
Thus we see that
\begin{equation}\label{munot}
\bmu_k(v_1\mcom v_n)=\s{(n-k)\mu_k}\hmu_k(v_1\mcom v_n).
\end{equation}
Using the notation of
section \ref{sect 5},
denote the restriction of $\hmu_k$ to $V^{k+l-1}$
by $\hmu_{kl}$. Set $n=k+l-1$. Denote
$\bmu_{kl}=\eta_l^{-1}\circ\delta_{kl}\circ\eta_{n}$, so that
$\bmu_{kl}$ is the restriction of $\bmu_k$ to $V^{n}$.
Then we  can express equation (\ref{munot}) in the form
$\bmu_{kl}=\s{(l-1)\mu_k}\mu_{kl}$.

More generally, if $\hdl$ is an arbitrary derivation on $T(W)$,
induced by the maps $\delta_k:V^k\ra V$, then it determines maps
$\mu_k:V^k\ra V$ , and $\bmu:T(V)\ra T(V)$, in a similar manner, and
we have
\begin{multline}
\bmu(v_1\mcom v_n)=
\sum
\begin{Sb}
1\le k\le n\\
\\
0\le i\le n-k
\end{Sb}
\s{ (v_1\mplus v_i)\mu_k+i(k-1) +(n-k)\mu_k}\\\times
v_1\mtns v_i\tns \mu_k(v_{i+1}\mcom v_{i+k})\tns v_{i+k+1}\mtns v_{n}.
\end{multline}
The condition that $\hdl$ is a codifferential on $T(W)$ can be expressed
in the form
\begin{equation}
\sum_{k+l=n+1}\delta_l\circ\delta_{kl}=0
\end{equation}
for all $n\ge1$.
This condition is equivalent to the condition
\begin{equation}
\sum_{k+l=n+1}\s{(l-1)\mu_k}\mu_l\circ\mu_{kl}=
\sum_{k+l=n+1}\mu_l\circ\bmu_{kl}
=0.
\end{equation}
We can express this condition in the form
\begin{multline}
\sum
\begin{Sb}
k+l=n+1
\\
0\le i\le n-k
\end{Sb}
\s{ (v_1\mplus v_i)\mu_k+i(k-1) +(n-k)\mu_k}\\\times
\mu_l(v_1\mcom v_i,\mu_k(v_{i+1}\mcom v_{i+k}),v_{i+k+1}\mcom v_{n})=0.
\end{multline}
When $\hdl$ is an odd codifferential, $\e{\mu_k}=k$,so the sign
in the expression above is simply
$(v_1\mplus v_i)k+i(k-1)+nk-k$. In \cite{ls,lm}, an element
$\mu\in\CV$ satisfying the equation
\begin{multline}
\sum
\begin{Sb}
k+l=n+1
\\
0\le i\le n-k
\end{Sb}
\s{ (v_1\mplus v_i)k+i(k-1) +nk-k}\\\times
\mu_l(v_1\mcom v_i,\mu_k(v_{i+1}\mcom v_{i+k}),v_{i+k+1}\mcom v_{n})=0.
\end{multline}
is called a
strongly homotopy associative algebra,
or \ainf\ algebra. We see that an \ainf\ algebra
structure on $V$ is nothing more than
an odd codifferential on $T(W).$ This
can also be expressed in terms of the bracket on $\CW$. An
odd element $\delta\in\CW$ satisfying $[\delta,\delta]=0$ determines a
codifferential $\hdl$ on $T(W)$.
More precisely, the condition
$[\delta,\delta]=0$ is equivalent to
the condition $\hdl^2=0$. Thus $\mu\in\CV$ determines an
\ainf\ algebra structure on $V$ when
$\delta=\eta\circ\mu\circ\eta^{-1}$ satisfies $[\delta,\delta]=0$.
This is not the same condition as $[\mu,\mu]=0$, nor even
the condition $\hmu^2=0$, although we shall have more to say about
this later.

Next set $V=\Pi U$, and $v_i=\pi u_i$.
Define $d_k:U^k\ra U$ by $d_k=\eta^{-1}\circ \mu_k\circ\eta$, and
$\brd:T(U)\ra T(U)$ by $\brd=\eta^{-1}\circ\mu\circ\eta$. Then
$\e{d_k}=\e{\mu_k}+(k-1)=\e{\delta_k}$.
Then by the same reasoning as above, we see that
\begin{multline}
d(u_1\mcom u_n)=
\sum_{i=0}^{n-k}
\s{(u_1\mplus u_i)d_k+n-nk}\\\times
u_1\mtns u_i\tns d_k(u_{i+1}\mcom u_{i+k})\tns u_{i+k+1}\mtns u_{n}.
\end{multline}
Arbitrary coderivations are treated in the same manner.
Continuing this process, let us suppose that $U=\Pi X$, and
$u_i=\pi x_i$. Define $m_k:X^k\ra X$ by
$m_k=\eta^{-1}\circ d_k\circ\eta$, and
$\bm:T(X)\ra T(X)$ by $\bm=\eta^{-1}\circ \brd \circ\eta$. Then
$\e{m_k}=\e{d_k}+(k-1)=\e{\mu_k}$.
Then we obtain that
\begin{multline}
\bm(x_1\mcom x_n)=
\sum_{i=0}^{n-k}
\s{(x_1\mplus x_i)m_k+i(k-1)+n-nk+(n-k)m_k}\\\times
x_1\mtns x_i\tns m_k(x_{i+1}\mcom x_{i+k})\tns x_{i+k+1}\mtns x_{n}.
\end{multline}
If we extend this process to an arbitrary coderivation, and consider
the signs which would result
from assuming that $\delta$ is an odd codifferential,
then we would obtain that
\begin{multline}\label{mysigns}
\sum
\begin{Sb}
0\le 1\le n-k\\
\\
k+l=n+1
\end{Sb}
\s{(x_1\mplus x_i)k+i(k-1)+n-k}\\\times
m_l(x_1\mcom x_i,m_k(x_{i+1}\mcom x_{i+k}),x_{i+k+1}\mcom x_{n})=0.
\end{multline}
The signs in the expression above agree with
the signs in the definition of a \ainf\ algebra
as given in \cite{ps2,kon}.

Finally, suppose that $X=\Pi Y$ and $y_i=\pi x_i$, and that we define
$\delta'_k=
\eta^{-1}\circ m_k\circ\eta$, and
$\bdl'=\eta^{-1}\circ \bm \circ\eta$. Then
$\e{\delta'_k}=\e{m_k}+(k-1)=\e{\delta_k}$.
Then we obtain that
\begin{multline}
\delta'(y_1\mcom y_n)=
\sum_{i=0}^{n-k}
\s{(y_1\mplus y_i)\delta_k}\\\times
y_1\mtns y_i\tns \delta'_k(y_{i+1}\mcom y_{i+k})\tns y_{i+k+1}\mtns y_{n}.
\end{multline}
Note that the signs occuring in this last expression are precisely the
same as for the original $\delta$. Thus it takes four rounds of parity
reversion to obtain the original signs. From this construction, it is
clear that the signs which arise in equation (\ref{mysigns}) are those
which would be obtained if we take $V=\Pi W$, and map $\CW$
into $\CV$ by setting $m=\eta\circ d\circ\eta^{-1}$ for
$d\in\CW$. Then a codifferential $d$ determines an \ainf\
structure on $V$ in the sense of \cite{ps2,kon}.

{}From these observations, we see that the two sign conventions for a
homotopy associative algebra originate because there are two natural
choices for the relationship between the space $W$, which carries
the structure of an odd differential with respect to the usual \zt-grading,
and $V$, which is its dual. One may choose either $W=\Pi V$, to get
the signs in \cite{ls,lm}, or $V=\Pi W$, to get the signs in
\cite{ps2,kon}. (Actually, one can vary the definition of $\eta$
to obtain both sets of signs from either one of these models.)

A natural question is why do we consider codifferentials
on the tensor coalgebra $T(W)$ of the parity reversion $W$
of $V$, rather than codifferentials on $T(V)$ in
the definition of an \ainf\ algebra? In fact, \ainf\ algebras are
generalizations of associative algebras, and an associative algebra
structure on $V$ is determined by a \ztz-graded odd codifferential
on $T(V)$. As a matter of fact, a map $m_2:V^2\ra V$ is an associative
multiplication exactly when $\hm_2$ is an odd codifferential.
The answer is that odd codifferentials with respect to the \zt-grading
have better properties with respect to the Lie bracket structure.

Let us examine the bracket structure on the space
$\coder(T(V))$.
In \cite{gers}, M. Gerstenhaber
defined a bracket on the space of cochains
of an associative algebra, which we shall
call the Gerstenhaber bracket.
When $V$ is concentrated in degree zero, in other words, in the non
\zt-graded case, the Gerstenhaber bracket is just the bracket of
coderivations, with the \Z-grading. Thus the Gerstenhaber bracket is
given by
\begin{equation}
[\ph_k,\psi_l]=\ph_k\psi_l-\s{(k-1)(l-1)}\psi_l\ph_k,
\end{equation}
for $\ph\in C^k(V)$ and $\psi\in C^l(V)$.
One of the main results of \cite{gers} is that the differential $D$ of
cochains in the cohomology of an associative algebra can be expressed
in terms of the bracket. It was shown that $D(\ph)=[\ph,m]$ where $m$
is the cochain representing the
associative multiplication. This formulation leads to a simple
proof that $D^2=0$, following from the properties of \Z-graded Lie
algebra. The associativity of $m$ is simply the condition $[m,m]=0$.
Let us recall the proof that an odd homogeneous element $m$ of a graded
Lie algebra satisfying $[m,m]=0$ gives rise to a differential $D$ on the
algebra by defining $D(\ph)=[\ph,m]$. In other words, we need to show
that $[[\ph,m],m]=0$. Recall that $m$ is odd when
$\ip{\e m}{\e m}=1$. The graded Jacobi bracket gives
\begin{equation}
[[\ph,m]m]
=[\ph[m,m]]+\s{\ip mm}[[\ph,m]m]=-[[\ph,m],m],
\end{equation}
which shows the desired result, since we are in characteristic zero.
Moreover, we point out that the Jacobi identity also shows that
$D([\ph,\psi])=[\ph,D(\psi)]+\s{\psi D}[D(\ph),\psi]$, so the differential
in the cohomology of an associative algebra acts as a graded derivation
of the Lie algebra, equipping $\CV$ with the structure of
a differential graded Lie algebra.

We wish to generalize the Gerstenhaber bracket to the \ainf\ algebra
case, where we are considering a more general codifferential $m$ on $T(V)$,
in such a manner that the bracket with $m$ yields a differential
graded Lie algebra structure on $\CV$. If we consider the
bracket of coderivations, then a problem arises when
the codifferential is not homogeneous. First of all,
$\hm^2=0$ does imply that $[m,m]=0$, but the converse is
not true in general. Secondly, if we define $D(\ph)=[\ph,m]$, then
we do not obtain in general that $D^2=0$. To see this, note that
$[m,m]=0$ is equivalent to $\sum_{k+l=n+1}[m_k,m_l]=0$ for all $n\ge 1$.
Let $\ph_p\in\hom(V^p,V)$. Then
\begin{multline}
[[\ph_p,m],m]_{n+p-1}=
\sum_{k+l=n+1}[[\ph_p,m_k],m_l]=\\
\sum_{k+l=n+1}[\ph_p,[m_k,m_l]]+
\s{\ip{m_k}{m_l}}
[[\ph_p,m_l],m_k]=\\
\sum_{k+l=n+1}\s{k+l+1}
[[\ph_p,m_l],m_k]=
\s{n}[[\ph_p,m],m]_{n+p-1}
\end{multline}
Thus we only obtain cancellation of terms when $n$ is odd. However, this
is sufficient to show that in the particular case where $m_k=0$ for all
even or all odd $k$, then $D^2=0$. In this case we also can show that
$\hm^2=0$ is equivalent to $[m,m]=0$ as well.

These problems occur because \ztz\ is not a good grading group.
Since these problems do not arise if we are considering \zt-graded
codifferentials (or \Z-graded codifferentials, in the \Z-graded case),
it is natural to consider a codifferential
on the parity reversion(or suspension)
of $V$, and using its properties to get a better
behaved structure.

For the remainder of this section, let us assume for definiteness
that $W=\Pi V$ and that $d\in\CW$ satisfies $\hd^2=0$.
In terms of the associative product induced on $\CW$, this
is the same condition as $d^2=0$. Because $W$ is \zt-graded, $d^2=0$
is equivalent to $[d,d]=0$ for $d$ odd, and moreover, $\CW$ is a differential
graded Lie algebra with differential $D(\ph)=[\ph,d]$. Thus we have no
problems with the differential on the $W$ side. Suppose that
$m_k=\eta_1\inv\circ d_k\circ\eta_k$ and
$\mu_l=\eta_1\inv\circ\delta_l\circ\eta_l$. Define a new bracket
$\brf$
on $\CV$ by
$\br{m_k}{\mu_l}=\eta_1\inv\circ[d_k,\delta_l]\circ\eta_{k+l-1}$. Then
it follows easily that
\begin{equation}
\br{m_k}{\mu_l}=
\s{(k-1)\mu_l}[m_k,\mu_l].
\end{equation}
Of course this new bracket no longer satisfies the graded antisymmetry
or graded Jacobi identity with respect to the inner product
we have been using
on \ztz. However, the bracket does satisfy these properties
with respect to the a different inner product on \ztz. We state
this result in the form of a lemma.
\begin{lma}\label{lma1}
Let $V$ be equipped with a \ztz-graded Lie bracket $[\cdot,\cdot]$
with respect to the inner product
$\ip{(\bar m,n)}{(\bar m',n')}=\bar m\bar m'+nn'$ on \ztz.
Then the bracket $\brf$ on $V$ given by
$\br uv=\s{\deg(u)|v|}[u,v]$
defines the structure of a \ztz-graded Lie algebra on $V$ with
respect to the inner product
$\ip{(\bar m,n)}{(\bar m',n')}=(\bar m+n)(\bar m'+n')$ on \ztz.

The bracket $\brtf$ given by
$\brt uv=\s{\deg(u)(|v|+\deg(v))}[u,v]$ also
defines the structure of a \ztz-graded Lie algebra on $V$ with
respect to the second inner product.
\end{lma}
The proof of the above lemma is straightforward, and will be omitted.

Because the bracket $\brf$ on $\CV$ really coincides with
the bracket of coderivations on $\CW$, we know that if
$[d,d]=0$, then $\br mm=0$, and we can define a differential on
$\CV$ by $D(\ph)=\br\ph m$, so that we can define a
cohomology theory for an \ainf\ algebra. We have been considering
the picture $W=\Pi V$. This means we have been considering \ainf\
algebras as defined in \cite{ls,lm}. If we consider instead the
picture $V=\Pi W$, then the bracket induced on $\CV$
by that on $\CW$ will be given by
$\br {m_k,\mu_l}=\s{(k-1)(\mu_l+l-1)}[m_k,\mu_l]$, which is the
second modified bracket described in the lemma. Thus we have
similar results. We shall call either one of these two brackets
the modified Gerstenhaber bracket.
The Hochschild cohomology of \ainf\ algebras was defined in
\cite{ps2} as the cohomology given by $D(\ph)=\br \ph m$, and
it was shown that this cohomology classifies the infinitesimal
deformations of an \ainf\ algebra. We shall not go into the details
here. It is important to note however, that unlike the cohomology
theory for an associative algebra, there is only one cohomology
group $H(V)$. The reason is that the image of $\ph\in C^n(V)$ under
the coboundary operator has a part in all $C^k(V)$ with $k\ge n$.
Only in the case of an associative or differential graded associative
algebra do we get a stratification of the cohomology.

Now let us suppose that $V$ is equipped with an inner
product $\ipf$.  The inner product induces an isomorphism between
$\hom(V^k,V)$ and $\hom(V^{k+1},\k)$, given by $\ph\mapsto\tilde\ph$,
where
\begin{equation}
\tilde\ph(v_1\mcom v_{k+1})=\ip{\ph(v_1\mcom v_k)}{v_{k+1}}.
\end{equation}
An element $\ph\in\hom(V^k,V)$ is said to be cyclic with respect to
the inner product if
\begin{equation}
\ip{\ph(v_1\mcom v_k)}{v_{k+1}}=
\s{k+v_1\ph}
\ip{v_1}{\ph(v_2\mcom v_{k+1})}.
\end{equation}
Then $\ph$ is cyclic if and only if $\tilde\ph$ is cyclic in the
sense that
\begin{equation}
\tilde\ph(v_1\mcom v_{k+1})=
\s{n + v_{k+1}(v_1\mplus v_k)}
\tilde\ph(v_{k+1},v_1\mcom v_k).
\end{equation}
The only differences between the definitions here and the definitions
given in section \ref{sect 2} is that $\tilde\ph$ is not completely
antisymmetric, and the isomorphism between $\hom(V^{k+1},\k)$ is to
the module $C^k(V)$, not just the submodule consisting of
cyclic elements.
If $m\in\CV$, then we say that $m$ is cyclic if $m_k$ is
cyclic for all $k$. If $m$ determines an \ainf\ algebra structure
on $V$, then we say that the inner product is invariant if $m$ is
cyclic with respect to the inner product.
The following lemma shows how to construct cyclic elements from
arbitrary elements of $\hom(\V^{n+1},\k)$.
\begin{lma}\label{lma2}
Suppose that $\tf\in\hom(V^{n+1},\k)$. Define
$C(\tf):V^{n+1}\ra\k$ by
\begin{multline}
C(f)(v_1\mcom v_{n+1})=\\
\sum_{0\le i\le n}
\s{(v_1\mplus v_i)(v_{i+1}\mplus v_{n+1})+ni}
f(v_{i+1}\mcom v_{i})
\end{multline}
Then $C(\tf)$ is cyclic. Furthermore, $C(\tf)=(n+1)\tf$ if $\tf$ is cyclic.
\end{lma}

We shall also need another lemma, which simplifies computations with
cyclic elements.
\begin{lma}\label{lma3}
Suppose that $\tilde f\in\hom(V^{n+1},\k)$ is cyclic,
$\alpha=v_1\mtns v_i$ and $\beta=v_i\mtns v_{n+1}$. Then
\begin{equation}
\tilde f(\alpha\tns\beta)=\s{\alpha\beta+in}\tilde f(\beta\tns\alpha)
\end{equation}

\end{lma}

Denote the submodule of cyclic elements in $\CV$ by $CC(V)$. The following
lemma records the fact that $CC(V)$ is a
graded Lie subalgebra of $\CV$, with
respect to the bracket of coderivations.

\begin{lma}\label{lma4}
Let $\ph_k\in C^k(V)$ and $\psi_l\in C^l(V)$. If $\ph$ and $\psi$ are
cyclic, then so is $[\ph,\psi]$. Moreover if $n=k+l-1$, then
\begin{multline}\label{cycrel}
\widetilde{[\ph,\psi]}(v_1\mcom v_{n+1})=\\
\sum
\begin{Sb}
0\le i\le n
\end{Sb}
\s{(v_1\mplus v_i)(v_{i+1}\mplus v_{n+1})+in}
\tilde\ph_k(\psi_l(v_{i+1}\mcom v_{i+l}),v_{i+l+1}\mcom v_i),
\end{multline}
where in the expression above, indices should be interpreted $\mod n+1$.
\end{lma}
\begin{pf}
The proof is straightforward, but involves a couple of tricks from
mod 2 addition in the signs, so we present it. Denote
$\rho=[\ph,\psi]$. Then
\begin{multline}
\tilde\rho(v_1\mcom v_{n+1})=\\
\sum_{0\le i\le k-1}
\s{(v_1\mplus v_i)\psi+i(l-1)}
\tilde\ph(v_1\mcom v_i,\psi(v_{i+1}\mcom v_{i+l}),v_{i+l+1}\mcom v_{n+1})\\
-\s{\ph\psi+(k-1)(l-1)}\times\\
\sum_{0\le i\le l-1} \s{(v_1\mplus v_i)\ph+i(k-1)}
\tilde\psi(v_1\mcom v_i,\ph(v_{i+1}\mcom v_{i+k}),v_{i+k+1}\mcom v_{n+1})
\end{multline}
Let us express the first term above, using the cyclicity of $\ph$.
\begin{multline}
\sum_{0\le i\le k-1}\s{r_i}
\tilde\ph(v_1\mcom v_i,\psi(v_{i+1}\mcom v_{i+l}),v_{i+l+1}\mcom v_{n+1})=\\
\sum_{0\le i\le k-1}\s{r_i+s_i}
\tilde\ph(\psi(v_{i+1}\mcom v_{i+l}),\mcom v_{i+l+1}\mcom v_i),
\end{multline}
where
\begin{eqnarray}
r_i&=&(v_1\mplus v_i)\psi+i(l-1)\\
s_i&=&(v_1\mplus v_i)(\psi + v_{i+1}\mplus v_{n+1}) +ki,
\end{eqnarray}
so that the coefficient is $r_i+s_i=(v_1\mplus v_i)(v_{i+1}\mplus v_{n+1})
+ni$.

The second term requires some additional manipulations.
\begin{multline}
\sum_{0\le i\le l-1} \s{r_i}
\tilde\psi(v_1\mcom v_i,\ph(v_{i+1}\mcom v_{i+k}),v_{i+k+1}\mcom v_{n+1})
=\\
\sum_{0\le i\le l-1} \s{r_i+s_i}
\tilde\psi(v_{i+k+1}\mcom v_i,\ph(v_{i+1}\mcom v_{i+k}))
=\\
\sum_{0\le i\le l-1} \s{r_i+s_i}
\ip{\psi(v_{i+k+1}\mcom v_i)}{\ph(v_{i+1}\mcom v_{i+k})}
=\\
\sum_{0\le i\le l-1} \s{r_i+s_i+t_i}
\ip{\ph(v_{i+1}\mcom v_{i+k})}{\psi(v_{i+k+1}\mcom v_i)}
=\\
\sum_{0\le i\le l-1} \s{r_i+s_i+t_i}
\tilde\ph(v_{i+1}\mcom v_{i+k},\psi(v_{i+k+1}\mcom v_i))
=\\
\sum_{0\le i\le l-1} \s{r_i+s_i+t_i+u_i}
\tilde\ph(\psi(v_{i+k+1}\mcom v_i)),v_{i+1}\mcom v_{i+k}),
\end{multline}
where
\begin{eqnarray}
r_i&=&\ph\psi+(k-1)(l-1)+1+(v_1\mplus v_i)\ph +i(k-1)\\
s_i&=&(v_1\mplus v_{i+k}+\ph)(v_{i+k+1}\mplus v_{n+1}+(i+1)l\\
t_i&=&(\psi+v_{i+k+1}\mplus v_i)(\ph+v_{i+1}\mplus v_{i+k})\\
u_i&=&(v_{i+1}\mplus v_{i+k})(\psi+v_{i+k+1}\mplus v_i)+k
\end{eqnarray}
so that
the sum of these coefficients is
$(v_i\mplus v_{i+k})(v_{i+k+1}\mplus v_{n+1})+in+kl$.
But now we use the fact that
$kn=k(l+k-1)=kl+k(k-1)=kl\mod 2$, so that $in+kl=(i+k)n$,
and re-indexing with $i=i+k$ allows us to express the second
term as
\begin{equation}
\sum_{k\le i\le n}
\s{(v_1\mplus v_i)(v_{i+1}\mplus v_{n+1})+ni}
\tilde\ph(\psi(v_{i+1}\mcom v_{i+l}),v_{i+l+1}\mcom v_i).
\end{equation}
Thus we have shown that equation (\ref{cycrel}) holds, and this is
enough to show that $\tilde\rho$ is cyclic, by lemma \ref{lma2}.
\end{pf}

Since $\br{\ph}{\psi}=\s{(k-1)\psi_l}[\ph,\psi]$, it follows that
the modified Gerstenhaber bracket of cyclic elements is also
cyclic.
Thus finally, we can state the main theorem, which allows us to define
cyclic cohomology of an \ainf\ algebra.

\begin{thm}\label{thm3}
i)Suppose that $V$ is a \zt-graded \k-module with an
inner product $\ipf$. Suppose that $\ph, \psi\in\CV$ are
cyclic. Then $\br{\ph}{\psi}$ is cyclic. Furthermore, the formula
below holds.
\begin{multline}
\widetilde{\br{\ph}{\psi}}(v_1\mcom v_{n+1})=
\sum
\begin{Sb}
k+l=n+1\\
\\
0\le i\le n
\end{Sb}
\s{(v_1\mplus v_i)(v_{i+1}\mplus v_{n+1})+in+(k-1)\psi_l}\times\\
\tilde\ph_k(\psi_l(v_{i+1}\mcom v_{i+l}),v_{i+l+1}\mcom v_i),
\end{multline}
where in the expression above, indices should be interpreted $\mod n+1$.
Thus the inner product induces a structure of a \ztz-graded
Lie algebra in the module $CC(V)$ consisting of all cyclic elements
in $\CV)$, by defining
$\br{\tilde\ph}{\tilde\psi}=\widetilde{\br{\ph}{\psi}}$.

ii) If $m$ is an \ainf\ structure on
$V$,
then there is a differential in $CC(V)$, given by
\begin{multline}
D(\tilde\ph)(v_1\mcom v_{n+1})=
\sum
\begin{Sb}
k+l=n+1\\
\\
0\le i\le n
\end{Sb}
\s{(v_1\mplus v_i)(v_{i+1}\mplus v_{n+1})+in +(k-1)l}\times\\
\tilde\ph_k(m_l(v_{i+1}\mcom v_{i+l}),v_{i+l+1}\mcom v_i),
\end{multline}
where in the expression above, $i$ should be interpreted $\mod n+1$.

iii) If the inner product is invariant, then
$D(\tilde\ph)=\br{\tilde\ph}{\tilde m}$. Thus $CC(V)$ inherits the
structure of a differential graded Lie algebra.
\end{thm}
\begin{pf}
The first statement follows from lemma \ref{lma4}, and the third
follows from the first two. The second assertion follows immediately
from the first when the inner product is invariant. The general case
is a routine verification, which we omit.
\end{pf}
We define $HC(V)$ to be the cohomology associated to the coboundary
operator on $CC(V)$. As in the case of Hochschild cohomology, there
is no stratification of the cohomology as $HC^n(V)$. Of course, $HC(V)$
is \zt-graded.
\section{Cohomology of \linf\ algbras}\label{sect 7}
There is a natural isomorphism $\eta$ between
$\bigwedge V$ and $\bigodot(\Pi V)$
which is given by
\begin{equation}
\eta(v_1\mwedge v_n)=
\s{(n-1)v_1\mplus v_{n-1}}
 \pi v_1\modot \pi v_n,
\end{equation}
Note that $\eta$ is neither
even nor odd. The restriction $\eta_k$ of $\eta$ to $\bigwedge^k V$ has parity
$k$. Of course, $\eta$ does preserve the exterior degree. For simplicity
in the following, let $W=\Pi V$ and let $w_i=\pi v_i$, and denote
$C(W)=\hom(\bigodot W,W)$, and $C(V)=\hom(\bigwedge V,V)$. We will
use notational conventions as in section \ref{sect 6}, so that
for $d\in C(W)$, $d_k$ will denote the restriction of this map to
$\bigwedge^k W$, $d_{lk}$ will denote the restriction of the
associated coderivation $\hd_l$ to $V^{k+l-1}$ etc.

The following lemma will be useful later on.
\begin{lma}\label{lma5}
Suppose that $\sigma$ is
a permutation of $n$ elements. Then
\begin{multline}
\s{(n-1)v_1\mplus v_{n-1}}\s{\sigma}\epsilon(\sigma;v_1\mcom v_n)=\\
\s{(n-1)v_{\sigma(1)}\mplus v_{\sigma(n-1)}}
\epsilon(\sigma;w_1\mcom w_n).
\end{multline}
\end{lma}
\begin{pf}
{}From the properties of the graded exterior algebra, we have
\begin{multline}
\eta(v_{\sigma(1)}\mwedge v_{\sigma(n)})=\\
\s{\sigma}\epsilon(\sigma;v_1\mcom v_n)
\eta(v_1\mcom v_n)=\\
\s{\sigma}\epsilon(\sigma;v_1\mcom v_n)
\s{(n-1)v_1\mplus v_{n-1}}
w_1\modot w_n.
\end{multline}
On the other hand, by direct substitution, we have
\begin{multline}
\eta(v_{\sigma(1)}\mwedge v_{\sigma(n)})=\\
\s{(n-1)v_{\sigma(1)}\mplus v_{\sigma(n-1)}}
w_{\sigma(1)}\modot w_{\sigma(n)}=\\
\s{(n-1)v_{\sigma(1)}\mplus v_{\sigma(n-1)}}
\epsilon(\sigma;w_1\mcom w_n)w_1\modot w_n
\end{multline}
Comparing the coefficients of the two expressions yields the desired
result.
\end{pf}
Let $d\in C(W)$ and  define $l_k=\eta_1\inv\circ d_k\circ\eta_k$.
so that
\begin{equation}
d_k(w_{\sigma(1)}\mcom w_{\sigma(k)}=
\s{(k-1)v_{\sigma(1)}\mplus v_{\sigma(k-1)}}
\pi m_k(v_{\sigma(1)}\mcom v_{\sigma(k-1)}).
\end{equation}
Let us abbreviate $\epsilon(\sigma;v_1\mcom v_n)$ by $\epsilon(\sigma;v)$
and $\epsilon(\sigma;w_1\mcom w_n)$ by $\epsilon(\sigma;w)$.
Let $n=k+l-1$.
Define $\bar l_{lk}=\eta_k\inv\circ d_{lk}\circ\eta_{n}$,
where $d_{lk}$ is given by
\begin{equation}
d_{lk}(w_1\mcom w_n)=
\sum_{\sigma\in\sh k{n-k}}
\epsilon(\sigma,w)
d_l(w_{\sigma(1)}\mcom w_{\sigma(l)})\odot w_{\sigma(l+1)}\modot w_{\sigma(n)}.
\end{equation}
We wish to compute $\bar l_{lk}$ in terms of $l_l$.
Now
\begin{multline}
l_{lk}(v_1\mcom v_n)=\\
\sum_{\sigma\in\sh l{n-l}}
\!\!\!
\epsilon(\sigma,w)
\s{(n-1)v_1\mplus v_{n-1}}
\eta_k\inv(
d_l(w_{\sigma(1)}\mcom w_{\sigma(l)})\odot w_{\sigma(l+1)}\modot w_{\sigma(n)})
\\=
\sum_{\sigma\in\sh l{n-l}}
\s{\sigma}\epsilon(\sigma,v)
\s{r}
\eta_k\inv(
d_l(w_{\sigma(1)}\mcom w_{\sigma(l)})\odot w_{\sigma(l+1)}\modot w_{\sigma(n)})
\\=
\sum_{\sigma\in\sh l{n-l}}
\s{\sigma}\epsilon(\sigma,v)
\s{r+s}
\eta_k\inv(
l_l(v_{\sigma(1)}\mcom v_{\sigma(l)})\odot w_{\sigma(l+1)}\modot w_{\sigma(n)})
\\=
\sum_{\sigma\in\sh l{n-l}}
\s{\sigma}\epsilon(\sigma,v)
\s{r+s+t}
l_l(v_{\sigma(1)}\mcom v_{\sigma(l)})\odot v_{\sigma(l+1)}\modot v_{\sigma(n)}
\\=
\sum_{\sigma\in\sh l{n-l}}
\s{\sigma}\epsilon(\sigma,v)
\s{(k-1)l_l}
l_l(v_{\sigma(1)}\mcom v_{\sigma(l)})\odot v_{\sigma(l+1)}\modot v_{\sigma(n)}.
\end{multline}

where
\begin{eqnarray}
r&=&(n-1)v_{\sigma(1)}\mplus v_{\sigma(n-1)}\\
s&=&(l-1)v_{\sigma(1)}\mplus v_{\sigma(l-1)}\\
t&=&(k-1)(l_l+v_{\sigma(1)}\mplus v_{\sigma(l-1)}
(k-2)v_{\sigma(l+1)}\mplus v_{\sigma(n-1)}
\end{eqnarray}

Thus we deduce immediately that $\bar l_{lk}=\s{(k-1)l_l})l_{lk}$, where
$l_{lk}$ is the restriction of the coderivation $\hat l_k$ to
$V^{k+l-1}$. This formula is identical to the formula we
deduced in section \ref{sect 6} connecting $\bm_{lk}$ to $\hm_{lk}$.

Now suppose that $d$ is an odd codifferential, so that
$d^2=0$.
This is equivalent to
$\sum_{k+l=n+1} d_k\circ d_{lk}=0$, which is equivalent to the
relations
$\sum_{k+l=n+1} \s{(k-1)l}l_k\circ l_{lk}=0$, since $\e{l_l}=l$.
This last relation can be put in the form
\begin{multline}
\sum
\begin{Sb}
k+l=n+1\\
\sigma\in\sh l{n-l}
\end{Sb}
\s{\sigma}\epsilon(\sigma,v)
\s{(k-1)l}
l_k(l_l(v_{\sigma(1)}\mcom v_{\sigma(l)}),
v_{\sigma(l+1)}\mcom v_{\sigma(n)}))=0.
\end{multline}
We say that the maps $l_k$ induce the structure of an \linf\ algebra,
or strongly homotopy Lie algebra
on $V$.
In \cite{ls}, the sign $k(l-1)$ instead of $(k-1)l$ appears in the
definition, but since $k(l-1)-(k-1)l=n+1$, this makes no difference
in the relations.

If we define $[\ph,\psi]$ to be the bracket of coderivations, then
we can define the modified Gerstenhaber bracket
$\br{\ph}{\psi}=\s{\deg\ph\e\psi}[\ph,psi]$. The definition of an
\linf\ algebra can be recast in terms of the bracket. In this language,
$l\in\CV$ determines an \linf\ structure on $V$ when $\br ll=0$.
The cohomology of an \linf\ algebra is defined to be the cohomology
on $\CW$ induced by $l$, in other words, $D(\ph)=[\ph, l]$.
This
definition makes $\CV$ a differential graded Lie algebra, with
respect to the second inner product on \ztz. These results are
completely parallel to the \ainf\ case.

In \cite{ps2}, the relationship between infinitesimal deformations
of an \ainf\ algebra and the cohomology of the \ainf\ algebra was
explored. The basic result is that the cohomology classifies the
infinitesimal deformations. Since we did not explore this matter
here, we shall discuss the parallel result for \linf\ algebras.

An infinitesimal deformation $l_t$ of an \linf\
algebra is given by taking
$l_t=l+t\lambda$, where the parity of $t$ is chosen so that
$(l_t)_k$ has parity $k$, so that we must have $\e{\lambda_k}=\e{t}+k$.
Since $t$ must have fixed parity, this determines the parity of
$\lambda_k$. The situation is more transparent if we switch to the
$W$ picture, so suppose that $l=\eta\inv\circ d\circ \eta$, and
$\lambda=\eta\inv\circ\delta\circ\eta$. Let $d_t=d +t\delta$, and
let us suppose that $d^2=0$, which is equivalent to $l$ giving
an \linf\ structure on $V$. Then $d_t$ is an infinitesimal
deformation of $d$ if $d_t^2=0$. Since $t$ is an infinitesimal
parameter, $t^2=0$, but we also want the parity of $d_t$ to be
odd, so that $\e t=1- \e\delta$. (We assume here that $\delta$
is homogeneous.) Now $d_t^2=0$ is equivalent to
$d^2+t\delta d+ dt\delta=0$. Also, $dt=\s{td}td$, the condition
reduces to $\delta d+\s{td}d\delta=0$. Now using the fact that
$\s{td}=\s{(1-\delta)}$, we see that $d_t$ is an infinitesimal
deformation precisely when $[\delta, d]=0$, in other words,
infinitesimal deformations are given by cocycles in the
cohomology of the \linf\ algebra. Trivial deformations are
more complicated, because we need to consider $d$ as a codifferential
of $T(W)$, and we define two codifferentials to be equivalent if
there is an automorphism of $T(W)$ which takes one of them to the
other. So trivial deformations depend on the structure of $T(W)$,
and are not simply given by taking an isomorphism of $W$ to itself.
However, one can show that $d_t$ is a trivial infinitesimal
deformation precisely when $\delta$ is a coboundary. Thus
the cohomology of $\CW$ classifies the infinitesimal deformations
of the \linf\ algebra. When we transfer this back to the $V$
picture, note that $m+\s{t}t\lambda$, is the deformed product
associated to $d_t$. But nevertheless, the condition for $l_t$
to be an \linf\ algebra still is that $\br \lambda\delta=0$.
We state
these results in the theorem below.
\begin{thm}
Let $l$ be an \linf\ algebra structure on $V$. Then the cohomology $H(V)$
of $\CV$ classifies the infinitesimal deformations of the \linf\
algebra.
\end{thm}

Moreover, let us suppose that $l$ is a Lie algebra structure on $V$.
Then the Lie algebra coboundary operator on $V$
coincides up to a sign with the \linf\ algebra coboundary operator
on $V$. This gives a nice interpretation of the cohomology of a
Lie algebra.
\begin{thm}
Let $l$ be an Lie algebra structure on $V$. Then the
Lie algebra cohomology $H(V)$
of $V$ classifies the infinitesimal deformations of the Lie algebra
into an \linf\ algebra.
\end{thm}
coboundary operator
Now let us address the case when $V$ is equipped with an inner product.
The bracket of coderivations  is the same bracket that was introduced
in equation \ref{nbra}
of section \ref{sect 2}, so that theorem \ref{th2} applies, and we
know that the bracket of two cyclic elements in $\CV$ is again cyclic.
Thus the modified Gerstenhaber bracket is cyclic as well. Thus we have
an analog of theorem \ref{thm3} for \linf\ algebras.
\begin{thm}
i)Suppose that $V$ is a \zt-graded \k-module with an
inner product $\ipf$. Suppose that $\ph, \psi\in\CV$ are
cyclic. Then $\br{\ph}{\psi}$ is cyclic. Furthermore, the formula
below holds.
\begin{multline}
\widetilde{\br{\ph}{\psi}}(v_1\mcom v_{n+1})=\\
\sum
\begin{Sb}
k+l=n+1\\
\\
\sigma\in\sh(l,k)
\end{Sb}
\s{\sigma}\epsilon(\sigma)\s{(k-1)\psi_l}
\tilde\ph_k(\psi_l(v_{\sigma(1)}\mcom v_{\sigma(l)}),
v_{\sigma(l+1)}\mcom v_{\sigma(n+1)}),
\end{multline}
Thus the inner product induces a structure of a \ztz-graded
Lie algebra in the module $CC(V)$ consisting of all cyclic elements
in $\CV)$, by defining
$\br{\tilde\ph}{\tilde\psi}=\widetilde{\br{\ph}{\psi}}$.

ii) If $l$ is an \linf\ structure on
$V$,
then there is a differential in $CC(V)$, given by
\begin{multline}
D(\tilde\ph)(v_1\mcom v_{n+1})=\\
\sum
\begin{Sb}
k+l=n+1\\
\\
\sigma\in\sh(l,k)
\end{Sb}
\s{\sigma}\epsilon(\sigma)\s{(k-1)l}
\tilde\ph_k(l_l(v_{\sigma(1)}\mcom v_{\sigma(l)}),
v_{\sigma(l+1)}\mcom v_{\sigma(n+1)}),
\end{multline}

iii) If the inner product is invariant, then
$D(\tilde\ph)=\br{\tilde \ph}{\tilde l}$. Thus $CC(V)$ inherits the
structure of a differential graded Lie algebra.
\end{thm}

We denote the cohomology given by the cyclic coboundary operator as
$HC(V)$.
Suppose that $V$ is an \linf\ algebra with an invariant inner product.
Then an infinitesimal deformation $l_t=l+t\ph$ preserves the inner product,
that is the inner product remains invariant under $l_t$, precisely when
$\ph$ is cyclic. Thus we see that cyclic cocycles correspond to infinitesimal
deformations of the \linf\ structure which preserve the inner product.
In a similar manner as before, cyclic coboundaries correspond to trivial
deformations preserving the inner product. Thus we have the following
classification theorem.
\begin{thm}
Let $l$ be an \linf\ algebra structure on $V$, with an invariant
inner product. Then the cyclic cohomology $HC(V)$
classifies the infinitesimal deformations of the \linf\
algebra preserving the inner product.
\end{thm}

Finally, we have an interpretation of the cyclic cohomology of a Lie
algebra.
\begin{thm}
Let $l$ be an Lie algebra structure on $V$, with an invariant
inner product. Then the
Lie algebra cyclic cohomology $H(V)$
of $V$ classifies the infinitesimal deformations of the Lie algebra
into an \linf\ algebra preserving the inner product.
\end{thm}
The author would like to thank Albert Schwarz,
Dmitry Fuchs and James Stasheff for reading this article and providing
useful suggestions.
\bibliographystyle{amsplain}
\ifx\undefined\bysame
\newcommand{\bysame}{\leavevmode\hbox to3em{\hrulefill}\,}
\fi

\end{document}